\documentclass[UTF-8,reqno]{amsart}
\usepackage{enumerate}
\usepackage{amssymb,url,color, booktabs}
\usepackage{mathrsfs}

\usepackage[nobysame]{amsrefs}
\BibSpec{article}{%
+{}{\PrintAuthors} {author}
+{,}{ \textrm} {title}
+{.}{ \textit} {journal}
+{,}{ \textbf} {volume}
+{}{ \parenthesize} {date}
+{,}{ } {pages}
+{.}{ arXiv:} {eprint}
+{.}{} {transition}
}
\BibSpec{book}{%
+{}{\PrintAuthors} {author}
+{,}{ \textit} {title}
+{.}{ \textrm} {series} %+{,}{ \textrm} {series}
+{,}{ Vol.} {volume} 
+{.}{ } {publisher}
+{,}{ } {date}
+{.}{} {transition}
}
\usepackage{color}
\usepackage[colorlinks=true]{hyperref}
\hypersetup{
    linkcolor=blue,          % color of internal links
    citecolor=red,        % color of links to bibliography
    filecolor=blue,      % color of file links
    urlcolor=cyan  
}

\numberwithin{equation}{section}

\newcommand{\be}{\begin{eqnarray}}
\newcommand{\ee}{\end{eqnarray}}
\newcommand{\ce}{\begin{eqnarray*}}
\newcommand{\de}{\end{eqnarray*}}
\newtheorem{theorem}{Theorem}[section]
\newtheorem{lemma}[theorem]{Lemma}
\newtheorem{remark}[theorem]{Remark}
\newtheorem{definition}[theorem]{Definition}
\newtheorem{proposition}[theorem]{Proposition}
\newtheorem{Examples}[theorem]{Example}
\newtheorem{corollary}[theorem]{Corollary}

\def\eps{\varepsilon}

\def\e{\mathrm{e}}

\def\p{\partial}

\def\[{{\Big[}}
\def\]{{\Big]}}
\def\<{{\langle}}
\def\>{{\rangle}}
\def\({{\Big(}}
\def\){{\Big)}}

\def\bx{{\mathbf{x}}}

\def\dif{{\mathord{{\rm d}}}}

\def\min{{\mathord{{\rm min}}}}

\def\no{\nonumber}
\def\={&\!\!=\!\!&}

\def\bB{{\mathbf B}}
\def\bC{{\mathbf C}}

\def\cR{{\mathcal R}}

\def\mE{{\mathbb E}}

\def\mN{{\mathbb N}}

\def\mP{{\mathbb P}}

\def\mR{{\mathbb R}}
\def\mS{{\mathbb S}}

\def\bB{{\mathbf B}}

\def\1{{\mathbf{1}}}

\def\sA{{\mathscr A}}
\def\sB{{\mathscr B}}

\def\sD{{\mathscr D}}

\def\sG{{\mathscr G}}

\def\sI{{\mathscr I}}

\def\sL{{\mathscr L}}

\def\sS{{\mathscr S}}

\def\E{\mathbb E}

\def\geq{\geqslant}
\def\leq{\leqslant}
\def\ge{\geqslant}

\def\div{\mathord{{\rm div}}}

\def\eps{\varepsilon}

\def\e{\mathrm{e}}

\def\p{\partial}

\def\[{{\Big[}}
\def\]{{\Big]}}
\def\<{{\langle}}
\def\>{{\rangle}}
\def\({{\Big(}}
\def\){{\Big)}}

\def\bx{{\mathbf{x}}}

\def\dif{{\mathord{{\rm d}}}}

\def\min{{\mathord{{\rm min}}}}

\def\no{\nonumber}
\def\={&\!\!=\!\!&}
\def\bt{\begin{theorem}}
\def\et{\end{theorem}}
\def\bl{\begin{lemma}}
\def\el{\end{lemma}}
\def\br{\begin{remark}}
\def\er{\end{remark}}
\def\bx{\begin{Examples}}
\def\ex{\end{Examples}}
\def\bd{\begin{definition}}
\def\ed{\end{definition}}
\def\bp{\begin{proposition}}
\def\ep{\end{proposition}}
\def\bc{\begin{corollary}}
\def\ec{\end{corollary}}

\def\geq{\geqslant}
\def\leq{\leqslant}
\def\ge{\geqslant}

\def\div{\mathord{{\rm div}}}

 \def\R{\mathbb R}
 \def\R{\mathbb R}

\def\<{\langle} \def\>{\rangle}

\allowdisplaybreaks

\begin{document}

\title[H\"older regularity and gradient estimates]{H\"older regularity and gradient estimates for SDEs driven by cylindrical $\alpha$-stable processes}

%\date{\today}

\author{Zhen-Qing Chen, Zimo Hao and Xicheng Zhang}

\address{Zhen-Qing Chen:
Department of Mathematics, University of Washington, Seattle, WA 98195, USA\\
Email: zqchen@uw.edu
 }
 
\address{Xicheng Zhang:
School of Mathematics and Statistics, Wuhan University,
Wuhan, Hubei 430072, P.R.China\\
Email: XichengZhang@gmail.com
 }

\address{Zimo Hao:
School of Mathematics and Statistics, Wuhan University,
Wuhan, Hubei 430072, P.R.China\\
Email: zimohao@whu.edu.cn
 }


\thanks{
Research of Z.-Q. Chen is partially supported
by  Simons Foundation grant 520542 and a Victor Klee Faculty Fellowship at UW.
Research of X. Zhang is partially supported by NNSFC grant of China (No. 11731009)  and the DFG through the CRC 1283 
``Taming uncertainty and profiting from randomness and low regularity in analysis, stochastics and their applications''. }

\begin{abstract}
We establish H\"older regularity and gradient estimates for the transition semigroup of the solutions to the following SDE:
$$
{\rm d} X_t=\sigma (t, X_{t-}){\rm d}  Z_t+b (t, X_t){\rm d} t,\ \ X_0=x\in\mR^d,
$$
where $( Z_t)_{t\geq 0}$ is a $d$-dimensional cylindrical $\alpha$-stable process with $\alpha \in (0, 2)$,
$\sigma (t, x):{\mathbb R}_+\times{\mathbb R}^d\to{\mathbb R}^d\otimes{\mathbb R}^d$ is bounded measurable, 
uniformly nondegenerate and Lipschitz continuous in $x$ uniformly in $t$,
and $b (t, x):{\mathbb R}_+\times{\mathbb R}^d\to{\mathbb R}^d$ is bounded $\beta$-H\"older continuous  in $x$ uniformly in $t$ 
with $\beta\in[0,1]$ satisfying $\alpha+\beta>1$. 
Moreover, we also show the existence and regularity of the distributional density of $X (t, x)$. Our proof is based on Littlewood-Paley's theory.
\end{abstract}

\maketitle

\section{\bf Introduction}

Consider the following stochastic differential equation (SDE) in $\mR^d$:
\begin{align}\label{SDE00}
\dif X_{s,t}=\sigma (t, X_{s,t})\dif W_t+b (t, X_{s,t})\dif t,\quad  t\geq s
\quad \hbox{with }   X_{s,s}=x\in\mR^d,
\end{align}
where $W$ is a $d$-dimensional standard Brownian motion,
$b (t, x):\mR_+\times\mR^d\to\mR^d$ is a bounded measurable function, and 
$\sigma (t, x):\mR_+\times\mR^d\to\mR^d\otimes\mR^d$ is a $d\times d$ matrix-valued measurable function.
Suppose that $\sigma$ satisfies the following uniformly elliptic condition:  there is some $c_0\geq 1$ so that 
\begin{align}
c^{-1}_0|\xi|\leq|\sigma (t, x)\xi|\leq c_0 |\xi|
\quad \hbox{for every }  t\geq 0 \hbox{ and  } 
 x,\xi\in\mR^d.\tag{{\bf H$^\sigma$}}
\end{align}
Under  ({\bf H$^\sigma$}) and $\lim_{|x-y|\to 0}\sup_{t\geq 0}|\sigma (t, x)-\sigma (t, y)|=0$, 
it is well known that for each starting point $(s,x)\in\mR_+\times\mR^d$, 
SDE \eqref{SDE00} admits a unique weak solution $X_{s,t}(x)$ (cf. \cite{St-Va}).
If, in addition, $\sigma (t, x)$ is H\"older continuous in $x$ uniformly in $t$, then  (cf. \cite{Zh-Zh18})
 for each $\varphi\in C^2_b(\mR^d)$,  
 $$
P_{s,t}\varphi(x):=\mE \left[ \varphi(X_{s,t}(x)) \right]
  $$
is $C^2$ in $x$, and uniquely 
solves the following backward Kolmogorov equation:
$$
\Big(\p_s+\tfrac{1}{2}\sigma^{ik} (s, \cdot) \sigma^{jk} (s, \cdot) \p_i\p_i+b^i (s, \cdot) \p_i\Big)P_{s,t}\varphi=0,\ \ s\leq t. 
$$ 
Here and below we use the usual Einstein convention: if an index appears twice in a product, then it will be summed automatically.

\medskip

A natural question is what kind regularity does $P_{s, t} \varphi$ have for solutions $X_{s, t}$ to
 SDEs driven by cylindrical stable processes, that is, 
for SDE \eqref{SDE00} with Brownian motion $W$ replaced by
a $d$-dimensional cylindrical $\alpha$-stable process $Z=(Z^1_t,\cdots,Z^d_t)_{t\ge 0}$? 
Here $\alpha\in(0,2)$ and $Z^1,\cdots, Z^d$ are independent $1$-dimensional 
$\alpha$-stable process with L\'evy measure $\dif z/|z|^{1+\alpha}$.  
In this paper we shall consider the following SDE driven by cylindrical  $\alpha$-stable process $Z$:
\begin{align}\label{SDE0}
\dif X_{s,t}=\sigma (t, X_{s,t-})\dif  Z_t+b (t, X_{s,t})\dif t,\quad  t\geq s
\quad \hbox{with }  X_{s,s}=x\in\mR^d.
\end{align}
Notice that the (time-dependent) generator for the solution $X_{s, t}$ of this SDE is  
\begin{align}\label{SL}
\sL^{\sigma,b}_t u(x):=\sum_{i=1}^d\mbox{p.v.}\int_{\mR}\Big(u(x+\sigma^{\cdot i} (t, x)z_i)-u(x)\Big)\frac{\dif z_i}{|z_i|^{1+\alpha}}+b^i (t, x)\p_i u(x),
\end{align}
where p.v. stands for the Cauchy principle value.
It is well known that if $\sigma$ and $b$ are Lipschtiz continuous in $x$ and uniformly in $t$, then there is a unique strong solution $X_{s,t}(x)$ to SDE \eqref{SDE0}.
When $b=0$,  Bass and Chen  \cite{Ba-Ch06} established the weak well-posedness for SDE \eqref{SDE0} under the assumption that $\sigma (t, x)=\sigma (x)$
is bounded, continuous and non-degenerate for each $x\in \R^d$.
 When $\sigma$ is Lipschtiz continuous in $x$  uniformly in $t$ and satisfies {\bf (H$^\sigma$)}, 
$b$ is $\beta$-order H\"older continuous
with $\alpha+\beta>1$,  Chen, Zhang and Zhao \cite{Ch-Zh-Zh} recently showed the strong well-posedness for SDE \eqref{SDE0} (see also \cite{Ch-So-Zh}).
Very recently, assuming that  $\alpha\in(0,1)$,  $b=0$ and $\sigma (t, x)=\sigma(x)$ is Lipschitz continuous satisfying {\bf (H$^\sigma$)},  
Kulczycki, Ryznar and Sztonyk \cite{Ku-Ry-Sz} showed the following H\"older estimate for the semigroup $P_t$ associated with SDE \eqref{SDE0}:
for any $\gamma\in(0,\alpha)$
and $T>0$, there is a positive constant $C=C(d, \alpha, c_0,  \| \nabla \sigma \|_\infty, \gamma, T)$ 
so that for all $t\in (0, T]$, 
$$
\|P_t\varphi\|_{\bC^\gamma}\leq 
C\,  t^{-{\gamma}/{\alpha}}\|\varphi\|_\infty,
$$
where $P_t\varphi(x):=\mE \left[ \varphi(X_t (x)) \right]$ and $\bC^\gamma$ is the space of bounded $\gamma$-H\"older continuous functions on $\R^d$.
 It was not known untill now if the above H\"older estimate holds when  $\alpha\in[1,2)$, nor was gradient estimate for $P_t$ for any $\alpha \in (0, 2)$.
 These properties will be addressed  in this paper under a more general setting. 
 When $Z$ is the rotationally invariant $\alpha$-stable process and $b,\sigma$ are Lipschitz continuous, 
the gradient estimate for $P_t$ was obtained in \cite{Wa-Xu-Zh} 
using  subordination technique.

  \medskip

The aim of this paper is to establish the following regularity estimates   for 
the transition semigroup $\{P_{s,t}; t\geq s\}$ of  the unique solution $\{X_{s, t} (x); t\geq s\}$ to SDE \eqref{SDE0}:
$$
P_{s,t}\varphi(x):=\mE \left[ \varphi(X_{s,t}(x)) \right].
$$
Note that  $\{P_{s,t}; t\geq s\}$ is the time-inhomogenous semigroup generated by the nonlocal operator 
$\sL^{\sigma,b}_t$ of \eqref{SL}.

\bt\label{Main}
Let $\alpha\in(0,2)$ and $\beta\in[0,1]$ with $\alpha+\beta>1$. 
Assume {\bf (H$^\sigma$)}, $\|\nabla\sigma\|_\infty\leq c_1$ for some $c_1>0$, and one of the following conditions holds:
$$
\mbox{{\rm (i)}  $b=0 $  and  $\beta=1$;\quad  {\rm (ii)}  $\alpha\in({1}/{2},2)$ and 
 $\sup_{t\geq 0}\|b(t,\cdot)\|_{\bC^\beta}\leq c_2$.}
$$
Let $\gamma\in[0,\alpha+\alpha\wedge\beta)$ and 
$\eta\in(-((\alpha+\beta-1)\wedge1),\gamma]$.
For any $T>0$, there exists a constant
 $C=C(d, c_0, c_1, c_2, \alpha, \beta, \gamma, \eta,  T)>0$ 
such that for all $0\leq s<t\leq T$,
\begin{align}\label{DW1}
\|P_{s,t}\varphi\|_{\bB^{\gamma}_{\infty,\infty}}\leq C(t-s)^{ (\eta-\gamma)/{\alpha}}\|\varphi\|_{\bB^\eta_{\infty,\infty}},
\end{align}
where $\bB^\eta_{p,q}$ is the usual Besov space.
In particular, we have the following assertions:
\begin{enumerate}[{\bf (A)}]
\item Let $\varphi\in\cup_{\eta<(\alpha+\beta-1)\wedge 1}\bB^{-\eta}_{\infty,\infty}$.  For any  $0\leq s<t$,
$P_{s,t}\varphi\in \cap_{\gamma<\alpha+\alpha\wedge\beta}\bB^{\gamma}_{\infty,\infty}$ 
and  the following backward Kolmogorov equation is satisfied:
for all $x\in\mR^d$,
\begin{align}\label{DW01}
P_{t_0,t}\varphi(x)=P_{t_1,t}\varphi(x)+\int^{t_1}_{t_0}\sL^{\sigma,b}_sP_{s,t}\varphi(x)\dif s,\ \ 0\leq t_0<t_1<t.
\end{align}
\item For $\alpha\in( {1}/{2},2)$, the following gradient estimate holds: for $0\leq s<t\leq T$,
\begin{align}\label{DW11}
\|\nabla P_{s,t}\varphi\|_\infty\leq C(t-s)^{-{1}/{\alpha}}\|\varphi\|_\infty.
\end{align}
\item For each $s<t$, the random variable $X_{s,t}(x)$ admits a density $p_{s,t}(x,\cdot)$ with
\begin{align}\label{REG1}
p_{s,t}(x,\cdot)\in \cap_{\eta<(\alpha+\beta-1)\wedge 1}\bB^\eta_{1,1}.
\end{align}
\end{enumerate}
\et

 \medskip
 
We would like to emphasize that for 
SDE \eqref{SDE0}  driven by cylindrical $\alpha$-stable process $Z$, 
 since the support of the L\'evy measure of $Z$ is concentrated along the 
coordinate axes,  it seems to be  quite difficult 
to obtain quantitative regularity results as stated in Theorem \ref{Main} by using 
methods from \cite{Ku-Ry-Sz} and \cite{Wa-Xu-Zh}. 
A new approach is needed to study regularity theory for SDEs driven by cylindrical stable processes. 
  A key ingredient in our approach is the use of  Littlewood-Paley's theory.

\medskip

\begin{remark}\rm 
\begin{itemize}
\item[(i)]
The precise definition of the Besov space $\bB^s_{p,q}$ is given in Definition \ref{iBesov} below. It is known that for non-integer 
$s>0$, $\|f\|_{\bB^s_{\infty,\infty}}\asymp\|f\|_{\bC^s}$. Hence \eqref{DW1} in particular yields 
that, under condition either (i) or (ii) of Theorem \ref{Main},  
for every $\alpha \in (0, 2)$, $T>0$,  and non-integer $ \gamma\in (0, \alpha+\alpha\wedge\beta)$,
there is a constant 
 $C=C(d, c_0, c_1, c_2,\alpha, \beta, \gamma, T)$ 
so that for all $0\leq s<t\leq T$, 
$$
\|P_{s,t}\varphi\|_{\bC^{\gamma}}\leq C(t-s)^{-{\gamma}/{\alpha}}\|\varphi\|_{\infty} . 
$$
This   significantly extends the main results of \cite{Ku-Ry-Sz} where $\alpha$ is restricted to be in $(0, 1)$, $\sigma(t, x)=\sigma(x)$,  
$b(t, x)\equiv 0$ and $\gamma \in (0, \alpha)$. 

\medskip
\item[(ii)]  Note that  $\sL^{\sigma,b}_s u(x)$ is pointwisely well defined for any $u\in \bC^\gamma$ with $\gamma >\max \{\alpha, 1\}$. 
 Under either condition (i) or (ii) of Theorem \ref{Main}, 
$\alpha +\alpha \wedge \beta >1$. Thus $\sL^{\sigma,b}_s P_{s, t} \varphi$ in \eqref{DW01}
is pointwisely well defined once it is established that  $P_{s,t}\varphi\in \cap_{\gamma<\alpha+\alpha\wedge\beta}\bB^{\gamma}_{\infty,\infty}$.

\medskip
\item[(iii)]
In the condition (ii) of Theorem \ref{Main},
due to some moment estimate, 
 $\alpha$ is required to be greater than ${1}/{2}$ (see Remark \ref{Br47} below).
This restriction also appears in the Schauder estimate of nonlocal PDEs in \cite{Ch-Me-Pr}. 
For variable coefficient $\sigma$, since we can only improve the regularity to $\alpha+\alpha\wedge 1$ even without drift $b$, 
we have to require $\alpha> {1}/{2}$ for gradient estimate also. An open problem is to drop the restriction $\alpha> {1}/{2}$.

\medskip
\item[(iv)] We note that when $b$ and $\sigma$ are time-independent, by a purely probabilistic method, 
Debussche and Fournier \cite{De-Fo} obtained the existence 
and low regularity of the densities for  SDE \eqref{SDE0} under some H\"older assumptions on $\sigma$ and $b$.
Compared with their results, for $\alpha\in[1,2)$, in the time independent case, the property \eqref{REG1} 
is covered by \cite[Theorem 1.1]{De-Fo}. 
However, for $\alpha\in( {1}/{2},1)$, the regularity \eqref{REG1} is better than \cite[Theorem 1.1]{De-Fo}.
The method in \cite{De-Fo} does not seem to be applicable to the time-dependent case and more general SDEs driven by Poisson random measures.
Our method is mostly analytic combined by some probabilistic argument and has more flexibility.

\medskip\item[(v)]  We point out that resolvent H\"older regularity can be established under a much weaker assumption on $\sigma$,
at least in the time-independent and driftless situation. 
Suppose that $\sigma (x)$ is continuous and satisfies condition {\bf (H$^\sigma$)}.
Then by \cite{Ba-Ch06} for each $x\in \R^d$, there is a unique weak solution   to  $\dif X_t = \sigma (X_{t-}) \dif Z_t$ with $X_0=x$,
where $ Z$ is a cylindrical stable process on $\R^d$.  Denote the law of $X$ with $X_0=x$ by $\mP_x$ and its mathematical expectation by $\E_x$.
It follows from \cite[Proposition 2.1]{BC10} and a scaling argument, that there are constants
 $c_3\geq c_4>0$ 
 so that for every $x_0\in \R^d$ and $r>0$,  
 \begin{equation}\label{e:1.9}
 \E_x  \left[ \tau_{B(x_0, r) }\right]  \leq 
 c_3
 r^\alpha \quad \hbox{for  } x\in B(x_0, r) 
 \end{equation}
 and
 \begin{equation}\label{e:1.10}
  \E_x \left[\tau_{B(x_0, r) } \right] \geq 
  c_4
    r^\alpha \quad \hbox{for } x\in B(x_0, r/2),
 \end{equation}
  where $B(x_0, r):=\{x:|x-x_0|<r\}$ and $\tau_{B(x_0, r)}=\inf\{t>0: X_t \notin B(x_0, r)\}$.
 Assume further that $\sigma (x)$ is uniformly continuous on $\R^d$, we know from
 \cite[Theorem, 2.9]{BC10} that there is a constant $\gamma \in (0, 1)$ that depends on the modulo of continuity of $\sigma(x)$,  the constant $c_0$ 
 in {\bf (H$^\sigma$)},  $d$ and $\alpha$ so that every bounded function that is harmonic in $B(x_0, r)$ is $\gamma$-H\"older continuous
 on $B(x_0, r/2)$.   
 This together with \eqref{e:1.9}-\eqref{e:1.10} and the proof of  
\cite[Proposition 2.4]{CCK} yields that there is some constant $\gamma_1 \in (0, 1)$ depending on  the modulo of continuity of $\sigma(x)$,  the constant $c_0$ 
 in {\bf (H$^\sigma$)}, $d$ and $\alpha$ so that for every  bounded function $\varphi$ on $\R^d$ and $\lambda >0$, 
 $$
 | R_\lambda \varphi (x) -R_\lambda \varphi (y)| \leq C |x-y|^{\gamma_1} \| \varphi \|_\infty  \quad \hbox{for } |x-y|\leq 1.
 $$
 Here 
 $R_\lambda \varphi (x):=\E_x \int_0^\infty e^{-\lambda t} \varphi (X_t ) \dif t$ is the $\lambda$-resolvent   of $\varphi$.  
\end{itemize}
\end{remark}

\medskip

This paper is organized as follows: In Section 2, we introduce some basic estimates for later use.
In Section 3, we present the estimates  of  Littlewood-Paley's type for the heat kernel of nonlocal operators with constant coefficients.
In Section 4, we show the regularity estimates for general nonlocal equations by freezing coefficients. In Section 5, we apply the a priori estimate obtained in 
Section 4 to show our main results. We use $:=$ as a way of definition.
Throughout this paper we shall use the following conventions and notations:
\begin{itemize}
	\item The letter $C$ denotes a constant, whose value may change in different places.
	\item We use $A\lesssim B$  to denote $A\leq C B$ for some unimportant constant $C>0$.
	\item $\mN_0:=\mN\cup\{0\}$, $\mR_+:=[0,\infty)$, $a\vee b:=\max(a,b)$, $a\wedge b:=\min(a,b)$.
	\item $\nabla_x:=\p_x:=(\p_{x_1},\cdots,\p_{x_d})$, $\p_i:=\p_{x_i}:=\p/\p x_i$. 
	\item For $x \in \R^d$ and $r>0$, we denote $B(x_0, r):=\{x\in\mR^d: |x-x_0|<r\}$ and $B_r:=B(0,r)$.
	\item For $p\in[1,\infty]$, we use $\|\cdot\|_p$ to denote the usual norm in $L^p(\mR^d,\dif x)$.
\end{itemize}

\medskip

\section{\bf Prelimiaries}

Let $\phi:\mR^d\to\mR^d$ be a bounded Lipschtiz function with
\begin{align}\label{SG0}
\|\nabla^\ell\phi\|_\infty\leq\kappa\leq  {1}/{2},  \quad \ell=0,1.
\end{align}
We shall use the following definitions: for a $C^1$-function $f:\mR^d\to\mR$,
\begin{align}\label{KH2}
\left\{
\begin{aligned}
&\Gamma_\phi(x):=x+\phi(x), \quad \sD^\phi_0 f(x):=f(x+\phi(x))-f(x+\phi(0)),\\
&\sD^\phi_1 f(x):=f(x+\phi(x))-f(x+\phi(0))-(\phi(x)-\phi(0))\cdot\nabla f(x).
\end{aligned}
\right.
\end{align}
The following lemma is  elementary (cf. \cite[Lemma 2.1]{Zh13}).
\bl
Under \eqref{SG0}, it holds that for any $x,y\in\mR^d$,
\begin{align}\label{Ga}
\tfrac{1}{2}|x-y|\leq |\Gamma_\phi(x)-\Gamma_\phi(y)|\leq 2|x-y|.
\end{align}
Moreover, there is a constant $C=C(d)>0$ such that for all $x\in\mR^d$,
\begin{align}\label{KF3}
|\det(\nabla\Gamma_\phi(x))-1-\div\phi(x)|\leq C\|\nabla\phi\|_\infty^2\leq C\kappa^2.
\end{align}
\el

For $\beta>0$, let $\bC^{\beta}$ be the space of  $\beta$-order H\"older continuous functions on $\mR^d$ with norm
$$
\|f\|_{\bC^\beta}:=\|f\|_\infty+\cdots+\|\nabla^{[\beta]} f\|_\infty+[\nabla^{[\beta]}  f]_{\bC^{\beta-[\beta]}}<\infty,
$$
where $[\beta]$ denotes the largest integer not exceeding $\beta$, and $\nabla^j$ stands for the $j$-order gradient,
and for $\gamma\in[0,1]$,
$$
[f]_{\bC^\gamma}:=\sup_h\|f(\cdot+h)-f(\cdot)\|_\infty/|h|^{\gamma}.
$$ 
We use the following convention: 
By $\bC^0$ we denote the space of bounded measurable functions. For two functions $f,g:\mR^d\to\mR$, let 
$$
\<f,g\>:=\int_{\mR^d} f(x) g(x)\dif x,
$$
whenever it is well defined.
The following lemma will play a crucial role in the proof of Theorem \ref{Main}.
\bl
Assume that $\phi:\mR^d\to\mR^d$ satisfies \eqref{SG0}. For any $\theta\in[0,1]$, there exists a constant $C=C(d,\theta)>0$ such that
for any $f\in L^\infty(\mR^d)$ and $g\in L^1(\mR^d)$ with $\nabla g\in L^1(\mR^d)$,
\begin{align}\label{KF6}
|\<\sD^\phi_0 f, g\>|\leq C\kappa^\theta\|f\|_\infty \left(\mu_0 (|g|)+\mu_\theta(|\nabla g|)^\theta\mu_\theta(|g|)^{1-\theta}\right),
\end{align}
and for any $f\in L^\infty(\mR^d)$ and $g\in L^1(\mR^d)$ with $\nabla g, \nabla^2 g\in L^1(\mR^d)$,
\begin{align}\label{KF8}
\begin{split}
|\<\sD^\phi_1 f, g\>|\leq 
C\kappa^{1+\theta}\|f\|_{\bC^\theta}\Big(  \sum_{j=0}^1\mu_j(|\nabla^j g|)
+\mu_{1+\theta}(|\nabla^2 g|)^\theta\mu_{1+\theta}(|\nabla g|)^{1-\theta}\Big),
\end{split}
\end{align}
where $\mu_\theta(\dif x):=(|x|\wedge 1)^\theta\dif x$ and $\mu_\theta(f):=\int_{\mR^d}f(x)\mu_\theta(\dif x)$.
\el
\begin{proof}
We first assume that
\begin{align}\label{Phi0}
\phi(0)=0.
\end{align} 
 {\it Step 1.}
Under \eqref{Phi0}, by a change of variable, we have
\begin{align}
|\<\sD^\phi_0 f, g\>|&=\left|\int_{\mR^d} \Big(f(\Gamma_\phi(x))-f(x)\Big)g(x)\dif x\right|\no\\
&=\left|\int_{\mR^d}f(x) \Big(g(\Gamma^{-1}_\phi(x))\det(\nabla\Gamma^{-1}_\phi(x))-g(x)\Big)\dif x\right|\no\\
&\leq \|f\|_\infty\Bigg(\int_{\mR^d}|g(\Gamma^{-1}_\phi(x))|\ |\det(\nabla\Gamma^{-1}_\phi(x))-1|\dif x\no\\
&\qquad\qquad\qquad+\int_{\mR^d}\big|g(\Gamma^{-1}_\phi(x))-g(x)\big|\dif x\Bigg).\label{KF4}
\end{align}
Since 
\begin{align}\label{Ga1}
\det(\nabla\Gamma^{-1}_\phi(x))=(\det\nabla\Gamma_\phi)^{-1}\circ \Gamma^{-1}_\phi(x),
\end{align}
by \eqref{KF3} and \eqref{SG0} we have
\begin{align}\label{KF5}
&\int_{\mR^d}|g(\Gamma^{-1}_\phi(x))|\cdot|\det(\nabla\Gamma^{-1}_\phi(x))-1|\dif x\no\\
&\quad=\int_{\mR^d}|g(x)|\cdot|\det(\nabla\Gamma_\phi(x))-1|\dif x 
\lesssim\kappa\int_{\mR^d}|g(x)|\dif x=\kappa\mu_0(|g|).
\end{align}  
On the other hand, noting that
$$
\big|g(\Gamma_\phi(x))-g(x)\big|\leq |\phi(x)|\int^1_0|\nabla g|(x+s\phi(x))\dif s=|\phi(x)|\int^1_0|\nabla g|(\Gamma_{s\phi}(x))\dif s,
$$
and since $\phi(0)=0$,
$$
|\phi(x)|=|\phi(x)-\phi(0)|\leq(\|\nabla\phi\|_\infty |x|)\wedge\|\phi\|_\infty\leq \kappa(|x|\wedge 1),
$$
we have by the change of variable again,  for any $\theta\in[0,1]$,
\begin{align*}
&\int_{\mR^d}\big|g(\Gamma^{-1}_\phi(x))-g(x)\big|\dif x=\int_{\mR^d}\big|g(x)-g(\Gamma_\phi(x))\big|\det(\nabla\Gamma_\phi(x))\dif x\\
&\qquad=\int_{\mR^d}\big|g(x)-g(\Gamma_\phi(x))\big|^\theta\big|g(x)-g(\Gamma_\phi(x))\big|^{1-\theta}\det(\nabla\Gamma_\phi(x))\dif x\\
&\qquad\stackrel{\eqref{Ga}}{\lesssim}\int_{\mR^d}\left(|\phi(x)|\int^1_0|\nabla g|(\Gamma_{s\phi}(x))\dif s\right)^\theta\Big(|g(\Gamma_\phi(x))|+|g(x)|\Big)^{1-\theta}\dif x\\
&\qquad\leq\kappa^\theta\int_{\mR^d}\left(\int^1_0|\nabla g|(\Gamma_{s\phi}(x))\dif s\right)^\theta\Big(|g(\Gamma_\phi(x))|+|g(x)|\Big)^{1-\theta}(|x|\wedge 1)^\theta\dif x\\
&\qquad\leq\kappa^\theta\left(\int^1_0\mu_\theta\Big(|\nabla g|\circ\Gamma_{s\phi}\Big)\dif s\right)^\theta
\left(\mu_\theta\Big(|g|\circ\Gamma_\phi+|g|\Big)\right)^{1-\theta},
 \end{align*}
where the last step is due to H\"older's inequality with respect to $\mu_\theta$.
Moreover,
\begin{align*}
\mu_\theta\Big(|\nabla^j g|\circ\Gamma_{s\phi}\Big)&=\int_{\mR^d}|\nabla^j g(x)|\cdot
\Big(|\Gamma^{-1}_{s\phi}(x)|\wedge 1\Big)^\theta\det(\nabla \Gamma^{-1}_{s\phi}(x))\dif x
\stackrel{\eqref{Ga}}{\lesssim}\mu_\theta\Big(|\nabla^j g|\Big),
\end{align*}
where $j=0,1$ and $s\in[0,1]$.
Hence,
\begin{align}\label{KF7}
\begin{split}
\int_{\mR^d}\big|g(\Gamma^{-1}_\phi(x))-g(x)\big|\dif x
&\lesssim \kappa^\theta\mu_\theta(|\nabla g|)^\theta\mu_\theta(|g|)^{1-\theta},
\end{split}
\end{align}
which together with \eqref{KF4} and \eqref{KF5} yields the desired estimate \eqref{KF6}.
\medskip
\\
 {\it Step 2.}
Under \eqref{Phi0}, as above, we have
\begin{align*}
\<\sD^\phi_1  f, g\>&=\int_{\mR^d} \Big(f(\Gamma_\phi(x))-f(x)-\phi(x)\cdot\nabla f(x)\Big)g(x)\dif x=\<f, \bar\sD^\phi_1 g\>,
\end{align*}
where
\begin{align*}
\bar\sD^\phi_1 g(x)&:=g(\Gamma^{-1}_\phi(x))\det(\nabla\Gamma^{-1}_\phi(x))-g(x)+(g\cdot\div\phi)(x)+(\phi\cdot\nabla g)(x)\\
&=g(\Gamma^{-1}_\phi(x))\Big(\det(\nabla\Gamma^{-1}_\phi(x))(1+(\div\phi)(\Gamma^{-1}_\phi(x)))-1\Big)\\
&\quad-(g\cdot\div\phi)(\Gamma^{-1}_\phi(x))\det(\nabla\Gamma^{-1}_\phi(x))+(g\cdot\div\phi)(x)\\
&\quad+g(\Gamma^{-1}_\phi(x))-g(x)+\phi(x)\cdot\nabla g(x)\\
&=:\sG_1g(x)+\sG_2g(x)+\sG_3g(x).
\end{align*}
In particular,
$$
\<\sD^\phi_1  f, g\>=\< f, \sG_1g\>+\< f, \sG_2g\>+\< f, \sG_3g\>.
$$
Notice that by \eqref{Ga1} and \eqref{KF3}, 
\begin{align*}
&|\det(\nabla\Gamma^{-1}_\phi(x))(1+(\div\phi)(\Gamma^{-1}_\phi(x)))-1|\\
&=\frac{|\det(\nabla\Gamma_\phi)-1-\div\phi|\circ \Gamma^{-1}_\phi(x)}{\det\nabla\Gamma_\phi\circ \Gamma^{-1}_\phi(x)}\leq C\kappa^2.
\end{align*}
Hence,
\begin{align}\label{AW1}
\<f,\sG_1 g\>\lesssim\kappa^2\|f\|_\infty \int_{\mR^d}|g\circ \Gamma^{-1}_\phi(x)|\dif x\lesssim\kappa^2 \|f\|_\infty\int_{\mR^d}|g(x)|\dif x.
\end{align}
Moreover, by the change of variable again, we have
\begin{align}
|\<f,\sG_2 g\>|&=\left|\int_{\mR^d}g(x)\cdot\div\phi(x)\cdot (f(\Gamma_\phi(x))-f(x))\dif x\right|\no\\
&\leq\|\div\phi\|_\infty\|\phi\|_\infty^\theta\|f\|_{\bC^\theta}\int_{\mR^d}|g(x)|\dif x
\leq\kappa^{1+\theta}\|f\|_{\bC^\theta}\mu_0(|g|).\label{AW2}
\end{align}
For $\sG_3 g(x)$, due to $\Gamma^{-1}_\phi(x)=x-\phi\circ\Gamma^{-1}_\phi(x)$, we have
\begin{align*}
\sG_3 g(x)&=\phi(x)\cdot\nabla g(x)-\phi\circ\Gamma^{-1}_\phi(x)\cdot\int^1_0 \nabla g(x-s\phi\circ\Gamma^{-1}_\phi(x))\dif s\\
&=\Big(\phi(x)-\phi\circ\Gamma^{-1}_\phi(x)\Big)\cdot\int^1_0 \nabla g(x-s\phi\circ\Gamma^{-1}_\phi(x))\dif s\\
&\quad+\phi(x)\cdot\int^1_0 \Big(\nabla g(x)-\nabla g(x-s\phi\circ\Gamma^{-1}_\phi(x))\Big)\dif s\\
&=:\sG_{31} g(x)+\sG_{32} g(x).
\end{align*}
For $\sG_{31} g(x)$, since
$$
|\phi(x)-\phi\circ\Gamma^{-1}_\phi(x)|\leq \|\nabla\phi\|_\infty|\phi\circ\Gamma^{-1}_\phi(x)|\leq\kappa^2(|x|\wedge1),
$$
we have
\begin{align*}
|\<f,\sG_{31} g\>|&\lesssim\kappa^2\|f\|_\infty\int^1_0\!\!\!\int_{\mR^d}(|x|\wedge1)|\nabla g|(x-s\phi\circ\Gamma^{-1}_\phi(x))\dif x\dif s\\
&\lesssim\kappa^2\|f\|_\infty\int^1_0\!\!\!\int_{\mR^d}(|\Gamma_\phi(x)|\wedge1)|\nabla g|(\Gamma_\phi(x)-s\phi(x))\dif x\dif s\\
&\lesssim\kappa^2\|f\|_\infty\int^1_0\!\!\!\int_{\mR^d}(|x|\wedge1)|\nabla g|(\Gamma_{(1-s)\phi}(x))\dif x\dif s\lesssim \kappa^2\|f\|_\infty\mu_1(|\nabla g|).
\end{align*}
For $\sG_{32} g(x)$, it is similar to \eqref{KF7} that for any $\theta\in[0,1]$,
\begin{align}\label{AW4}
|\<f,\sG_{32} g\>|\lesssim\kappa^{1+\theta}\|f\|_\infty\mu_{1+\theta}(|\nabla^2 g|)^\theta\mu_{1+\theta}(|\nabla g|)^{1-\theta}.
\end{align}
Combining \eqref{AW1}-\eqref{AW4}, we obtain \eqref{KF8} under \eqref{Phi0}.
\medskip\\
 {\it Step 3.}
In the general case, without assuming \eqref{Phi0}, if we define
$$
\bar f(x):=f(x+\phi(0)),\ \bar\phi(x):=\phi(x)-\phi(0),
$$
then
$$
\sD^\phi_0 f(x)=\sD^{\bar\phi}_0 \bar f(x),\ \ \sD^\phi_1 f(x)=\sD^{\bar\phi}_1 \bar f(x)+\bar\phi(x)\cdot\nabla (\bar f-f)(x).
$$
For \eqref{KF6}, it follows by (i). For \eqref{KF8}, by (ii), it remains to make the following estimate:
\begin{align*}
|\<\bar\phi\cdot\nabla (\bar f-f),g\>|&=\left|\int_{\mR^d}\Big(\div \bar\phi(x) g(x)+\bar\phi(x)\cdot\nabla g(x)\Big)(\bar f-f)(x)\dif x\right|\\
&\leq\int_{\mR^d}\Big(\kappa |g(x)|+\kappa(|x|\wedge 1)|\nabla g(x)|\Big)|\phi(0)|^\theta\cdot[f]_{\bC^\theta}\dif x\\
&\leq\kappa^{1+\theta}\Big(\mu_0(|g|)+\mu_1(|\nabla g|)\Big)[f]_{\bC^\theta}.
\end{align*}
The proof is complete.
\end{proof}

Let $\sS(\mR^d)$ be the Schwartz space of all rapidly decreasing functions on $\mR^d$, and $\sS'(\mR^d)$ 
the dual space of $\sS(\mR^d)$ called Schwartz generalized function (or tempered distribution) space. Given $f\in\sS(\mR^d)$, 
its Fourier transform $\hat f$ and inverse Fourier transform $\check f$ are defined by
\begin{align*}
\hat f(\xi):=(2\pi)^{-d/2}\int_{\mR^d} \e^{-i\xi\cdot x}f(x)\dif x, \quad
\check f(x):=(2\pi)^{-d/2}\int_{\mR^d} \e^{i\xi\cdot x}f(\xi)\dif\xi.
\end{align*}
Let $\phi_0$ be a radial $C^\infty$-function on $\mR^d$ with 
$$
\phi_0(\xi)=1\ \mbox{ for } \ |\xi|\leq 1\ \mbox{ and }\ \phi_0(\xi)=0\ \mbox{ for } \ |\xi|\geq2.
$$
Define for $j\in\mN$,
$$
\phi_j(\xi):=\phi_0(2^{-j}\xi)-\phi_0(2^{1-j}\xi).
$$
It is easy to see that for $j\in\mN$, $\phi_j(\xi)=\phi_1(2^{1-j}\xi)\geq 0$ and
$$
{\rm supp}\phi_j\subset B_{2^{j+1}}\setminus B_{2^{j-1}},\ \  \sum_{j=0}^{k}\phi_j(\xi)=\phi_0(2^{-k}\xi)\to 1,\ \ k\to\infty.
$$
\bd\label{iBesov}
For given $j\in\mN_0$, the block operator $\cR_j$ is defined on $\sS'(\mR^d)$ by
\begin{align}\label{Def2}
\cR_jf(x):=(\phi_j\hat f)\check{\,\,}(x)=\check\phi_j* f(x)=2^{d (j-1)}\int_{\mR^d}\check\phi_1(2^{j-1}y) f(x-y)\dif y.
\end{align}
For any $s\in\mR$ and $p,q\in[1,\infty]$, the Besov space $\bB^s_{p,q}$ is defined by
$$
\bB^s_{p,q}:=\left\{f\in\sS'(\mR^d): \|f\|_{\bB^s_{p,q}}:=\left(\sum_{j\in\mN_0} 2^{sq j} \|\cR_j f\|^q_p \right)^{1/q}<\infty\right\}.
$$
\ed
\br\label{Re24}
It is well known that for $0<s\notin\mN$ (cf. \cite{Tr92}): 
\begin{align}\label{LJ1}
\|f\|_{\bB^s_{\infty,\infty}}\asymp\|f\|_{\bC^s}.
\end{align}
Moreover, let $\bar\cR_j:=\cR_{j-1}+\cR_j+\cR_{j+1}$ with convention $\cR_{-1}=0$. 
Since $\phi_{j-1}+\phi_j+\phi_{j+1}=1$ 
 on $B_{2^{j+1}}\setminus B_{2^{j-1}}$,
we have
$$
(\phi_{j-1}+\phi_j+\phi_{j+1})\phi_j=\phi_j.
$$
Consequently, 
\begin{align}\label{TR1}
\cR_j\bar\cR_j=\bar\cR_j\cR_j=\cR_j.
\end{align}
\er
The following commutator estimate is proven in \cite{Ch-Zh-Zh} (see also \cite{Ba-Ch-Da}).
\bl
Let $p,p_1,p_2\in[1,\infty]$ with $\frac{1}{p}=\frac{1}{p_1}+\frac{1}{p_2}$.
For any $\beta\in(0,1)$ and $\gamma\in(-\beta,0]$, there is a constant $C=C(d,\beta,\gamma)>0$ such that
for all $j\in\mN$,
\begin{align}
\|[\cR_j, f]g\|_p\leq C2^{-j(\beta+\gamma)}\|f\|_{\bB^{\beta}_{p_2,\infty}}\|g\|_{\bB^{\gamma}_{p_1,\infty}},\ j\in\mN_0,\label{GS1}
\end{align}
where $[\cR_j,f]g:=\cR_j(fg)-f\cR_j g$.
\el

We also need the following Gronwall inequality of Volterra type (cf. \cite[Lemma 2.2]{Zh10}).
\bl\label{Le26}
Let $f\in L^1_{loc}(\mR_+;\mR_+)$ and $T>0$. Suppose that for some $\gamma,\beta\in[0,1)$ and $C_1,C_2>0$,
$$
f(t)\leq C_1t^{-\beta}+C_2\int^t_0(t-s)^{-\gamma}f(s)\dif s,\ \ t\in(0,T].
$$
Then there is a constant $C_3=C_3(C_2,T,\gamma,\beta)>0$ such that
$$
f(t)\leq C_3 C_1 t^{-\beta},\ \ t\in(0,T].
$$
\el

\medskip

 \section{\bf Gradient estimates for heat kernel of nonlocal operator with constant coefficient}

Fix $\alpha\in(0,2)$. Let $( Z_t)_{t\geq 0}$ be an $\alpha$-stable process with L\'evy measure
\begin{align}\label{Le}
\nu(A)=\int^\infty_0\frac{\dif r}{r^{1+\alpha}}\int_{\mS^{d-1}}\1_A(r\omega)\pi(\dif\omega),\ A\in\sB(\mR^d\setminus\{0\}),
\end{align}
where $\pi$ is a finite measure over the unit sphere $\mS^{d-1}$. Note that the $\alpha$-stable process $Z$ has the scaling property
\begin{align}\label{DG1}
(\lambda^{-1/\alpha} Z_{\lambda t})_{t\geq 0}\stackrel{(d)}{=}( Z_t)_{t\geq 0},\ \ \forall\lambda>0,
\end{align}
and for any $\gamma>\alpha>\beta\geq 0$.
\begin{align}\label{DG2}
\int_{|z|<1}|z|^\gamma\nu(\dif z)+\int_{|z|>1}|z|^\beta\nu(\dif z)<\infty.
\end{align}
Let $N(\dif t, \dif z)$ be the associated Poisson random measure, that is,
$$
N([0,t]\times A):=\sum_{s\in(0,t]}\1_A(\Delta  Z_s),\ \ A\in\sB(\mR^d\setminus\{0\}), t>0,
$$
where $\Delta  Z_s:= Z_s- Z_{s-}$ is the jump size of $Z$ at time $s$.
Let $\phi(t,z):\mR_+\times\mR^d\to\mR^d$ be a measurable function with
$$
|\phi(t,z)|\leq C|z|, \ \ (t,z)\in\mR_+\times\mR^d.
$$
We consider the following time-inhomogenous L\'evy process: for $0\leq s\leq t<\infty$,
$$
X^\phi_{s,t}:=
\left\{
\begin{aligned}
&\int^t_s\!\!\int_{\mR^d}\phi(r,z)N(\dif r,\dif z),&\alpha\in(0,1),\\
&\int^t_s\!\!\int_{|z|\leq 1}\phi(r,z)\tilde N(\dif r,\dif z)+\int^t_s\!\!\int_{|z|>1}\phi(r,z)N(\dif r,\dif z),&\alpha=1,\\
&\int^t_s\!\!\int_{\mR^d}\phi(r,z)\tilde N(\dif r,\dif z),&\alpha\in(1,2),
\end{aligned}
\right.
$$
where $\tilde N(\dif r,\dif z):=N(\dif r,\dif z)-\dif r\nu(\dif z)$ is the compensated Poisson random measure. Define
\begin{align}\label{KG1}
P^\phi_{s,t}f(x):=\mE f(x+X^\phi_{s,t}),\ \ f\in C^2_b(\mR^d).
\end{align}
By It\^o's formula, one has
\begin{align}\label{HS1}
\p_t P^\phi_{s,t}f=\sL^{(\alpha)}_{\phi_t} P^\phi_{s,t}f,
\end{align}
where 
 $\phi_t (z):=\phi(t,z)$ 
and
$$
\sL^{(\alpha)}_{\phi_t} f(x):=\int_{\mR^d} \Big(f(x+\phi(t,z))-f(x)-\phi^{(\alpha)}(t,z)\cdot\nabla f(x)\Big)\nu(\dif z),
$$
and
\begin{align}\label{PH}
\phi^{(\alpha)}(t,z):=\Big(\1_{\alpha\in(1,2)}+\1_{\alpha=1}\1_{|z|\leq 1}\Big)\phi(t,z).
\end{align}

\bl(Duhamel's formula)\label{Le37}
Let $\varphi\in\bC^0$ and $f\in  L^\infty_{loc}(\mR_+;\bC^0)$. Define
$$
u(t,x):=P^\phi_{0,t}\varphi(x)+\int^t_0P^\phi_{s,t} f(s,x)\dif s.
$$
Then $u\in L^\infty_{loc}(\mR_+;\bC^0)$ uniquely solves $\p_t u=\sL^{(\alpha)}_{\phi_t} u+f$ in the distributional sense,
that is, for any $\psi\in C^2_c(\mR^d)$ and $t>0$,
\begin{align}\label{DD3}
\<u(t),\psi\>=\<\varphi,\psi\>+\int^t_0\<u(s),\sL^{(\alpha)}_{-\phi_s}\psi\>\dif s+\int^t_0\<f(s),\psi\>\dif s.
\end{align}
\el
\begin{proof}
Recall that $\bC^0$ stands for the space of bounded measurable functions. Clearly, $u\in L^\infty_{loc}(\mR_+;\bC^0)$.
Let $(\rho_\eps)_{\eps\in(0,1)}$ be a family of mollifiers in $\mR^d$ with support in $B_\eps$. Define
$$
u_\eps(t):=u(t)*\rho_\eps,\ \ \varphi_\eps:=\varphi*\rho_\eps,\ \ f_\eps(t):=f(t)*\rho_\eps.
$$
Clearly,
$$
u_\eps(t,x)=P^\phi_{0,t}\varphi_\eps(x)+\int^t_0P^\phi_{s,t} f_\eps(s,x)\dif s.
$$
By \eqref{HS1} and the integration by parts, one sees that
\begin{align}\label{DD2}
u_\eps(t,x)=\varphi_\eps(x)+\int^t_0\sL^{(\alpha)}_{\phi_s} u_\eps(s,x)\dif s+\int^t_0f_\eps(s,x)\dif s.
\end{align}
Hence, for any $\psi\in C^2_c(\mR^d)$, 
$$
\<u_\eps(t),\psi\>=\<\varphi_\eps,\psi\>+\int^t_0\<u_\eps(s),\sL^{(\alpha)}_{-\phi_s}\psi\>\dif s+\int^t_0\<f_\eps(s),\psi\>\dif s.
$$
Taking $\eps\to 0$, 
we obtain \eqref{DD3}. On the other hand, if we take
$\psi(\cdot)=\rho_\eps(x-\cdot)$ in \eqref{DD3}, then we get \eqref{DD2} and the uniqueness follows.
\end{proof}

Below we always make the following assumptions:
\begin{enumerate}[{\bf (H$^{\phi,\nu}$)}]
\item $\phi$ satisfies the following nondegeneracy conditions
$$
0<c^\phi_0:=\inf_{\omega_0\in\mS^{d-1}}\inf_{t,\lambda>0}\int_{\mS^{d-1}}\frac{|\omega_0\cdot\phi(t,\lambda \omega)|^2}{\lambda^2}\pi(\dif\omega),\ \
\sup_{t,z}\frac{|\phi(t,z)|}{|z|}=:c^\phi_1<\infty,
$$
and
\begin{align}\label{GH1}
\1_{\alpha=1}\int_{R_1<|z|<R_2}\phi(t,z)\nu(\dif z)=0,\ \ 0<R_1<R_2.
\end{align}

\end{enumerate}
Notice that  for all $0\leq s<t<\infty$,
\begin{align}\label{MA1}
c_0^\phi=c_0^{\phi_{s,t}}\mbox{ and }c_1^\phi=c_1^{\phi_{s,t}},
\end{align}
and by \eqref{DG1} and \eqref{GH1},
\begin{align}\label{Sc0}
X^\phi_{s,t}\stackrel{(d)}{=}(t-s)^{-1/\alpha}X^{\phi_{s,t}}_{0,1},
\end{align}
where 
$$
\phi_{s,t}(r,z):=(t-s)^{1/\alpha}\phi(s+r(t-s),(t-s)^{-1/\alpha}z).
$$

The following lemma can be proved as in \cite{Ch-Zh18b} (see also \cite{Ch-Me-Pr}).
For the readers' convenience, we provide a detailed proof here.
\bl\label{Le31}
Under {\bf (H$^{\phi,\nu}$)}, for each $0\leq s<t<\infty$, $X_{s,t}$ admits a $C^\infty$-smooth density 
$p^\phi_{s,t}(x)$ which satisfies that for any $\beta\in[0,\alpha)$ and $n\in\mN_0$,
\begin{align}\label{He1}
\int_{\mR^d}|x|^\beta|\nabla^n p^\phi_{s,t}(x)|\dif x\leq C(t-s)^{\frac{\beta-n}{\alpha}},
\end{align}
where $C=C(c_0^\phi,c_1^\phi,n,d,\alpha,\beta)>0$. Moreover, for each $0\leq s<t<\infty$ and $x\in\mR^d$,
\begin{align}\label{Sca}
p^\phi_{s,t}(x)=(t-s)^{-d/\alpha}p^{\phi_{s,t}}_{0,1}((t-s)^{-1/\alpha}x).
\end{align}
\el
\begin{proof}
First of all, by \eqref{Sc0} one sees that \eqref{Sca} holds.
Thus by \eqref{MA1}, 
it suffices to prove \eqref{He1} for $s=0$ and $t=1$. We only consider the case $\alpha\in(0,1)$ and write
$$
X^\phi_{0,1}=\int^1_0\!\!\int_{|z|<1}\phi(r,z)N(\dif r,\dif z)+\int^1_0\!\!\int_{|z|\geq 1}\phi(r,z)N(\dif r,\dif z)=:\widetilde X^\phi_{0,1}+\bar X^\phi_{0,1}.
$$
Note that the characteristic function of $\widetilde X^\phi_{0,1}$ is given by
$$
\mE\e^{{\rm i}\xi\cdot\widetilde X^\phi_{0,1}}=\e^{\psi_1(\xi)},
$$
where for $\delta\in(0,\infty]$,
$$
\psi_\delta(\xi):=\int^1_0\!\!\int_{|z|<\delta} \Big(\e^{\mathrm{i}\xi\cdot\phi(t,z)}-1\Big)\nu(\dif z)\dif t.
$$
We claim that there is a constant $c_2>0$ such that for all $\xi\in\mR^d$,
\begin{align}\label{star}
{\rm Re}\psi_1(\xi)\leq -c_2\Big(|\xi|^2\wedge|\xi|^\alpha\Big).
\end{align}
Indeed, by a change of variable, we have
\begin{align*}
{\rm Re}\psi_\infty(\xi)&=\int^1_0\!\!\int^\infty_0\frac{1}{r^{1+\alpha}}\int_{\mS^{d-1}} \Big(\cos(\xi\cdot\phi(t,r\omega))-1\Big)\pi(\dif \omega)\dif r\dif t\\
&=|\xi|^\alpha\int^1_0\!\!\int^\infty_0\frac{1}{r^{1+\alpha}}\int_{\mS^{d-1}} \Big(\cos(\xi\cdot\phi(t,r\omega/|\xi|))-1\Big)\pi(\dif \omega)\dif r\dif t\\
&\leq |\xi|^\alpha\int^1_0\!\!\int^\delta_0\frac{1}{r^{1+\alpha}}\int_{\mS^{d-1}} \Big(\cos(\xi\cdot\phi(t,r\omega/|\xi|))-1\Big)\pi(\dif \omega)\dif r\dif t.
\end{align*}
Note that
\begin{align}\label{GN1}
\lim_{r\to 0}\frac{1-\cos(r)}{r^2}=\frac{1}{2}.
\end{align}
By {\bf (H$^{\phi,\nu}$)}, since $|\xi\cdot\phi(t,r\omega/|\xi|)|\leq c^\phi_1 r$,
we can choose $\delta$ small enough so that
\begin{align*}
{\rm Re}\psi_\infty(\xi)&\leq-\frac{|\xi|^\alpha}{3}\int^1_0\!\!\int^\delta_0\frac{1}{r^{1+\alpha}}\int_{\mS^{d-1}}|\xi\cdot\phi(t,r\omega/|\xi|)|^2\pi(\dif \omega)\dif r\dif t\\
&\leq-c^\phi_0\frac{|\xi|^\alpha}{3}\int^\delta_0\frac{r^2}{r^{1+\alpha}}\dif r,
\end{align*}
and therefore, there are constants $c_3>0$ and $M>0$ such that for all $|\xi|>M$,
\begin{align*}
{\rm Re}\psi_1(\xi)&={\rm Re}\psi_\infty(\xi)+|\psi_\infty(\xi)-\psi_1(\xi)|\\
&\leq -c^\phi_0\frac{|\xi|^\alpha}{3}\int^\delta_0\frac{r^2}{r^{1+\alpha}}\dif r+2\int_{|z|>1}\nu(\dif z)\leq-c_3|\xi|^\alpha.
\end{align*}
On the other hand, by \eqref{GN1} and {\bf (H$^{\phi,\nu}$)}, for $\delta$ small enough, we also have for $|\xi|\leq M$,
\begin{align*}
{\rm Re}\psi_1(\xi)&\leq\int^1_0\!\!\int_{|z|\leq\delta}
\Big(\cos(\xi\cdot\phi(t,z))-1\Big)\nu(\dif z)\dif t\\
&\leq -c_4\int^1_0\!\!\int_{|z|\leq\delta}|\xi\cdot\phi(t,z)|^2\nu(\dif z)\dif t\\
&=-c_4|\xi|^2\int^1_0\!\!\int^\delta_0\frac{1}{r^{1+\alpha}}\int_{\mS^{d-1}} 
|\bar\xi\cdot\phi(t,r\omega)|^2\pi(\dif \omega)\dif r\dif t\\
&\leq-c_4c^\phi_0|\xi|^2\int^\delta_0\frac{r^2}{r^{1+\alpha}}\dif r,
\end{align*}
where $\bar\xi:=\xi/|\xi|$.
This proves the claim \eqref{star}.

By \eqref{star}, it is easy to see that $\e^{\psi_1(\xi)}\in\sS(\mR^d)$ is a Schwartz function. 
Thus the random variable $\widetilde X^\phi_{0,1}$ 
has a $C^\infty$-smooth density $\rho \in \sS(\mR^d)$. Noting that for $\beta<\alpha$,
$$
\mE|\bar X^\phi_{0,1}|^\beta\leq C(c^\phi_1,\beta,\alpha,d)<\infty,\ \ p^\phi_{0,1}(x)=\mE\rho\Big(x+\bar X^\phi_{0,1}\Big),
$$
we have
\begin{align*}
&\int_{\mR^d}|x|^\beta|\nabla^np^\phi_{0,1}(x)|\dif x
\leq\mE\int_{\mR^d}|x|^\beta|\nabla^n\rho|\big(x+\bar X^\phi_{0,1}\big)\dif x\\
&\quad\lesssim\int_{\mR^d}|x|^\beta|\nabla^n\rho|(x)\dif x+\mE|\bar X^\phi_{0,1}|^\beta\int_{\mR^d}|\nabla^n\rho|(x)\dif x<\infty.
\end{align*}
The desired estimate \eqref{He1} for $s=0$ and $t=1$ follows.
\end{proof}

 The following is a key lemma, which is similar to \cite[Lemma 3.1]{Ha-Wu-Zh}.
\bl\label{Le33}
Assume {\bf (H$^{\phi,\nu}$)} and let $p_{s,t}(x):=p^\phi_{s,t}(x)$ be as  in Lemma \ref{Le31}.
\begin{enumerate}[(i)]
\item For any $n\in\mN_0$, there is a constant $C>0$ such that for all $0\leq s<t<\infty$,
\begin{align}\label{KF10}
\|\nabla^n\cR_0 p_{s,t}\|_1\leq C.
\end{align}
\item For any $n\in\mN_0$, $\vartheta\geq 0$ and $\beta\in[0,\alpha)$, there is a constant $C>0$ such that  for all $0\leq s<t<\infty$ and $ j\in\mN$,
\begin{align}\label{KF2}
 m_\beta\Big(|\nabla^n\cR_j p_{s,t}|\Big)\leq C 2^{(n-\vartheta)j}(t-s)^{-\vartheta/\alpha}((t-s)^{\beta/\alpha}+2^{-j\beta}),
\end{align}
where $ m_\beta(\dif x):=|x|^\beta\dif x$ and $ m_\beta(\varphi):=\int_{\mR^d}\varphi(x) m_\beta(\dif x)$.
\item For any $n\in\mN_0$, $\beta,\gamma\in[0,\alpha)$,  there is a constant $C>0$ such that 
for all nonnegative measurable $f:\mR_+\to\mR_+$, $0\leq s<t<\infty$, $j\in\mN$,
\begin{align}\label{KF1}
\int^t_0 m_\beta(|\nabla^n\cR_j p_{s,t}|)f(s)\dif s\leq C 2^{(n-\gamma-\beta)j}\int^t_0(t-s)^{-\frac{\gamma}{\alpha}}f(s)\dif s.
\end{align}
\end{enumerate}
\el
\begin{proof}
(i) Let $n\in\mN_0$.
By the definition of $\cR_0$, we have
$$
\|\nabla^n\cR_0p_{s,t}\|_1\lesssim \|\cR_0 p_{s,t}\|_1\leq\|\check\phi_0\|_1\|p_{s,t}\|_1=\|\check\phi_0\|_1<\infty.
$$
(ii)
Since the support of $\phi_1$ is contained in the annulus, for any $k\in\mN_0$,
$$
\Delta^{-k}\check\phi_1:=(|\xi|^{-2k}\phi_1(\xi))\check{\,}\in\sS(\mR^d).
$$ 
Fix $j\in\mN$. For simplicity of notation, we write
$$
\tilde p_{0,1}(x):=p^{\phi_{s,t}}_{0,1}(x),\ \ \hbar:=(t-s)^{- {1}/{\alpha}}2^{-j}.
$$
Thus by \eqref{Sca} and the change of variable, for any $k\in\mN_0$, we have
\begin{align*}
\cR_j p_{s,t}(x)&=(t-s)^{-d/\alpha}\int_{\mR^d}\tilde p_{0,1}((t-s)^{-1/\alpha}2^{-j}y)\check\phi_1(2^{j}x-y)\dif y\\
&=\hbar^{d+k} 2^{jd}\int_{\mR^d}(\Delta^k\tilde p_{0,1})(\hbar y)(\Delta^{-k}\check\phi_1)(2^{j}x-y)\dif y.
\end{align*}
Therefore,  for any $k\in\mN_0$,
\begin{align*}
& m_\beta\Big(|\nabla^n\cR_j p_{s,t}|\Big)=\int_{\mR^d}\!\!|x|^\beta |\nabla^n\cR_j p_{s,t}(x)|\dif x\\
&=\!\hbar^{d+k} 2^{j(n-\beta)}\!\!\!\int_{\mR^d}\!\!|x|^\beta\!\left|\int_{\mR^d}\!\!(\Delta^k\tilde p_{0,1})(\hbar y)(\nabla^n\Delta^{-k}\check\phi_1)(x-y)\dif y\right|\!\dif x\\
&\lesssim \hbar^{d+k} 2^{j(n-\beta)}\int_{\mR^d}|x|^\beta|\Delta^k \tilde p_{0,1}(\hbar y)|\dif y\int_{\mR^d}|\nabla^n\Delta^{-k}\check\phi_1(x)|\dif x\\
&+\hbar^{d+k} 2^{j(n-\beta)}\int_{\mR^d}|\Delta^k\tilde p_{0,1}(\hbar y)|\dif y\int_{\mR^d}|x|^\beta|\nabla^n\Delta^{-k}\check\phi_1(x)|\dif x\\
&\stackrel{\eqref{He1}}{\lesssim} \hbar^{k} 2^{j(n-\beta)}\Big(\hbar^{-\beta}+1\Big),
\end{align*}
which in turn gives \eqref{KF2} by simple interpolation.
\medskip
\\
(iii) Let $\sI$ be the left hand side of \eqref{KF1}. Without loss of generality, we assume $t> 2^{-\alpha j}$. We make the following decomposition:
\begin{align*}
\sI=\left(\int^t_{(t-2^{-\alpha j})\vee 0}+\int^{(t-2^{-\alpha j})\vee 0}_0\right)\mu_\beta\Big(|\nabla^n\cR_j p_{s,t}|\Big) f(s)\dif s=:\sI_1+\sI_2.
\end{align*}
For $\sI_1$, by \eqref{KF2} with $\vartheta=\gamma$, we have 
\begin{align*}
\sI_1&\lesssim 2^{(n-\gamma)j}\int^t_{(t-2^{-\alpha j})\vee 0}(t-s)^{-\frac{\gamma}{\alpha}}\Big(2^{-j}+(t-s)^{\frac{1}{\alpha}}\Big)^\beta f(s)\dif s\\
&\leq 2^{(n-\gamma-\beta)j}\int^t_0(t-s)^{-\frac{\gamma}{\alpha}} f(s)\dif s.
\end{align*}
For $\sI_2$, by \eqref{KF2} with $\vartheta=\gamma+\beta$, we have 
\begin{align*}
\sI_2&\lesssim 2^{(n-\gamma-\beta)j} \int^{(t-2^{-\alpha j})\vee 0}_0(t-s)^{-\frac{\gamma+\beta}{\alpha}}\Big(2^{-j}+(t-s)^{\frac{1}{\alpha}}\Big)^\beta f(s)\dif s\\
&\leq 2^{(n-\gamma-\beta)j}\int^t_0(t-s)^{-\frac{\gamma}{\alpha}} f(s)\dif s.
\end{align*}
Combining the above two estimates, we obtain \eqref{KF1} for $j\in\mN$.
\end{proof}

\medskip

\section{\bf Regularity estimate for nonlocal equations}

In this section we fix $\alpha\in(0,2)$ and consider the following time-dependent nonlocal operator:
$$
\sL^{(\alpha)}_{\phi_t}u(x):=\int_{\mR^d} \Big(u(x+\phi(t,x,z))-u(x)-\phi^{(\alpha)}(t,x,z)\cdot\nabla u(x)\Big)\nu(\dif z),
$$
where $\nu$ takes the form \eqref{Le}, 
  $ \phi(t,x,z):\mR_+\times\mR^d\times\mR^d\to\mR^d$
is a measurable function, and 
 $$
\phi^{(\alpha)}(t, x, z):=\Big(\1_{\alpha\in(1,2)}+\1_{\alpha=1}\1_{|z|\leq 1}\Big)\phi(t,x, z).
$$
Recall that $\pi$ is the finite measure on $\mS^{d-1}$ associated with the L\'evy measure $\nu$ in \eqref{Le}.
We  assume for some $c_0,c_1>0$,
\begin{align}\label{GH13}
|\phi(t,x,z)|\leq c_0|z|,\quad |\nabla_x\phi(t,x,z)|\leq c_1|z|,
\end{align}
and
\begin{align}\label{GH12}
\inf_{\omega_0\in\mS^{d-1}}\inf_{x\in\mR^d}\inf_{t,\lambda>0}\int_{\mS^{d-1}}\frac{|\omega_0\cdot\phi(t,x,\lambda \omega)|^2}{\lambda^2}\pi(\dif\omega)=:c_2>0,
\end{align}
and
\begin{align}\label{GH11}
\1_{\{ \alpha=1 \}}\int_{R_1<|z|<R_2}\phi(t,x,z)\nu(\dif z)=0 \quad \hbox{for any } 0<R_1<R_2.
\end{align}
Clearly, $\sL^{(\alpha)}_{\phi_t}u(x)$ is  well defined pointwisely
if $u\in\bC^\gamma$ 
 for some $\gamma>\alpha$.  
Let $b:\mR_+\times\mR^d\to\mR^d$ be a measurable function and satisfy that for some $\beta\in[0,1]$,
\begin{align}\label{BB}
|b(t,x)|\leq c_3,\ \ |b(t,x)-b(t,y)|\leq c_3|x-y|^\beta.
\end{align}
The following parameter set will be used for stating the dependence of a constant.
$$
\Theta:=(d,\alpha,c_0,c_1,c_2,c_3,\beta).
$$
Fix $s\geq 0$. Consider the following nonlocal equation:
\begin{align}\label{PDE0}
\p_t u_s=\sL^{(\alpha)}_{\phi_t}u_s+\1_{\{\alpha>1/2\}} b\cdot\nabla u_s
\quad \text{for}\ t\geq s \quad  \hbox{with }   u_s(s,x)=\varphi(x).
\end{align}
In order to introduce the classical solution of \eqref{PDE0}, we define
$$
\bC^{\gamma-}:=\cap_{\gamma'<\gamma}\bC^{\gamma'}.
$$
\bd\label{Def41}
Fix $s\geq 0$ and $\gamma>\alpha\vee 1$. For $\varphi\in \bC^{\gamma-}$, we call a function 
$u_s\in C([s,\infty); \bC^{\gamma-})$
a classical solution of nonlocal-PDE \eqref{PDE0} with initial value $\varphi$ at time $s$ if for all $t\geq s$ and $x\in\mR^d$,
$$
u_s(t,x)=\varphi(x)+\int^t_s\Big(\sL^{(\alpha)}_{\phi_r}+\1_{\{\alpha>1/2\}} b\cdot\nabla\Big) u_s(r,x)\dif r.
$$
\ed

By the proof of \cite[Theorem 6.1]{CHXZ}, maximum principle holds for classical solutions of \eqref{PDE0}
and so classical solution  to \eqref{PDE0}   is unique. 
Fix $\gamma>\alpha\vee 1$ and suppose $u_s (t, x)$ is the classical solution to \eqref{PDE0}. 
To explicitly reflect its dependence on its initial value $\varphi$ at time $s$, 
we write  
$$
Q_{s,t}\varphi(x):=u_s(t,x).
$$
It follows from the uniqueness of classical solution to \eqref{PDE0} that 
for any $s<r<t$,
\begin{align}\label{GQ1}
Q_{s,t}\varphi=Q_{r,t}Q_{s,r}\varphi.
\end{align}

\medskip

We first establish the following a priori regularity estimate.

\bt\label{Th31}
Let $\alpha\in(0,2)$, $\beta\in[0,1]$ with $\alpha+\beta>1$ and $\gamma\in[0,\alpha+\alpha\wedge \beta)$.
Under conditions \eqref{GH13}-\eqref{BB},  for any $T>0$ and $\eta\in(-((\alpha+\beta-1)\wedge1),\gamma]$,
there is a constant $C=C(T, \Theta,\gamma, \eta)>0$ such that for any $0\leq s<t\leq T$, $\varphi\in\bC^2$ and any
classical solution $Q_{s,t}\varphi(x)=u_s(t,x)$ of the nonlocal-PDE \eqref{PDE0},
\begin{align}\label{DJ1}
\|Q_{s,t}\varphi\|_{\bB^\gamma_{\infty,\infty}}\leq C(t-s)^{\frac{\eta-\gamma}{\alpha}}\|\varphi\|_{\bB^\eta_{\infty,\infty}}.
\end{align}
\et
To show \eqref{DJ1},  
we use the freezing coefficient argument. Without loss of generality, we assume $s=0$ and write $u(t,x)=Q_{0,t}\varphi(x)$.
Fix $y\in\mR^d$ and let $\sS^y$ be the set of all solutions $\theta^y$ of the following ODE:
$$
\dot \theta^y(t)=-b(t,\theta^y(t)),\ \ \theta^y(0)=y.
$$
It is well known that for any $t>0$ (e.g. \cite[Lemma 6.5]{Ha-Wu-Zh}),
\begin{align}\label{AQ1}
\cup_{y\in\mR^d}\cup_{\theta^y\in\sS^y}\{\theta^y(t)\}=\mR^d.
\end{align}
Define
$$
u^y(t,x):=u(t,x+\theta^y(t)),\ \ \phi^y (t, x,z):=\phi^y(t,x,z):=\phi(t,x+\theta^y(t),z)
$$
and
$$
b^y(t,x):=b(t,x+\theta^y(t))-b(t,\theta^y(t)).
$$
It is easy to see that
$$
\p_t u^y+\sL^{(\alpha)}_{\phi^y_t}u^y+\1_{\{\alpha>1/2\}} b^y\cdot\nabla u^y=0,\ \ u^y(0,x)=\varphi^y(x):=\varphi(x+y).
$$
In the following, for notional simplicity
we drop the superscript $y$ from $u$, $\varphi, b$.  With this notation, $u$ satisfies 
$$
\p_t u+\sL^{(\alpha)}_{\phi_t}u+\1_{\{\alpha>1/2\}} b\cdot\nabla u=0,\ \ u(0)=\varphi,
$$
and $b$ satisfies
\begin{align}\label{KF9}
|b(t,x)|\leq 2c_3(|x|^\beta\wedge 1).
\end{align}
Next we freeze $\phi$ at point $0$. Define $\psi (t, z):=\phi (t, 0,z)$ and 
$$
\sA u:=\sL^{(\alpha)}_{\phi_t} u-\sL^{(\alpha)}_{\psi_t}u.
$$
Then we have
\begin{align}\label{tstar}
\p_t u+\sL^{(\alpha)}_{\psi_t} u+\sA u+\1_{\{\alpha>1/2\}} b\cdot\nabla u=0,\ \ u(0)=\varphi.
\end{align}
Let $P^{\psi}_{s,t}$ be defined by \eqref{KG1} in terms of $\phi=\psi$ and $p^{\psi}_{s,t}$ the corresponding heat kernel, that is,
$$
P^{\psi}_{s,t}f(x)=\int_{\mR^d}p^{\psi}_{s,t}(x-y)f(y)\dif y.
$$
Since $u$ is a classical solution of \eqref{tstar}, by   
Lemma \ref{Le37},  
\begin{align*}
u(t,x)=P^{\psi}_{0,t}\varphi^y(x)+\int^t_0P^{\psi}_{s,t}(\sA u)(s,x)\dif s+\1_{\{\alpha>1/2\}}\int^t_0P^{\psi}_{s,t}(b\cdot\nabla u)(s,x)\dif s.
\end{align*}
For $j\in\mN_0$,   
acting on both sides of the above equation by $\cR_j$, we obtain
\begin{align}\label{KL1}
\begin{split}
\cR_j u(t,0)&=\cR_jP^{\psi}_{0,t}\varphi^y(0)+\int^t_0\cR_j P^{\psi}_{s,t}(\sA u)(s,0)\dif s\\
&\quad+\1_{\{\alpha>1/2\}}\int^t_0\cR_jP^{\psi}_{s,t}(b\cdot\nabla u)(s,0)\dif s.
\end{split}
\end{align}

\bl\label{Le42}
For any $T>0$ and $\eta\leq\gamma$, there is a constant $C=C(T,\gamma,\eta,\Theta)>0$ 
such that  for all $j\in\mN_0$, $t\in(0,T]$ and $y\in\mR^d$,
$$
2^{\gamma j}|\cR_jP^{\psi}_{0,t}\varphi^y(0)|\leq Ct^{\frac{\eta-\gamma}{\alpha}}\|\varphi\|_{\bB^\eta_{\infty,\infty}},
$$
where $\varphi^y(x):=\varphi(x+y)$.
\el
\begin{proof}
By definition and H\"older's inequality, we have for any $\eta\leq\gamma$,
\begin{align*}
|\cR_jP^{\psi}_{0,t}\varphi^y(0)|&=\left|\int_{\mR^d}\cR_j p^{\psi}_{0,t}(-x)\varphi^y(x)\dif x\right|
\stackrel{\eqref{TR1}}{=}\left|\int_{\mR^d}\cR_j p^{\psi}_{0,t}(-x)\bar\cR_j\varphi^y(x)\dif x\right|\\
&\leq\|\cR_j p^{\psi}_{0,t}\|_1\|\bar\cR_j\varphi\|_\infty\lesssim 2^{-\gamma j}t^{-\frac{\gamma-\eta}{\alpha}}\|\varphi\|_{\bB^\eta_{\infty,\infty}},
\end{align*}
where the last step is due to Lemma \ref{Le33} and the definition of $\bB^\eta_{\infty,\infty}$.
\end{proof}

\bl\label{Le43}
For any $T>0$, $\gamma\in[0,\alpha)$ and $\eps\in(0,\alpha-\gamma)$, there is a constant $C=C(T,\eps,\gamma,\Theta)>0$ such that for all $j\in\mN_0$ and $t\in(0,T]$,
$$
2^{\gamma j}\int^t_0|\cR_jP^{\psi}_{s,t}\sA u(s,0)|\dif s\leq 
\int^t_0 (t-s)^{-\frac{\gamma+\eps}{\alpha}}\|u(s)\|_{\bB^{(\alpha-1)\vee 0+\eps}_{\infty,\infty}}\dif s.
$$
\el
\begin{proof}
We only prove the estimate for $\alpha\in(1,2)$. The case $\alpha\in(0,1]$ is similar and easier.
Since the time variable does not play any essential role in the following calculations, if there is no confusions, we shall drop the time variable  for simplicity of notation.
 Let $\delta>0$ be a fixed small number, which will be determined below. 
Since $\alpha\in(1,2)$, by definition we can make the following decomposition:
$$
\sA u(x)=\int_{|z|\leq\delta} \sD^{\phi(\cdot,z)}_1 u(x)\nu(\dif z)+\int_{|z|>\delta} \sD^{\phi(\cdot,z)}_1 u(x)\nu(\dif z)=:
\sA_\delta u(x)+\bar\sA_\delta u(x),
$$
where $\sD^{\phi(\cdot,z)}_1 u(x)$ is defined by (see \eqref{KH2})
$$
\sD^{\phi(\cdot,z)}_1 u(x):=u(x+\phi(x,z))-u(x+\phi(0,z))-(\phi(x,z)-\phi(0,z))\cdot\nabla u(x).
$$
We first treat $\sA_\delta u$.
Notice that  by definition and Fubini's theorem,
\begin{align*}
\cR_jP^{\psi}_{s,t}\sA_\delta u(s,0)&=\int_{\mR^d} \cR_j p^{\psi}_{s,t}(-x)\sA_\delta u(s,x)\dif x\\
&=\int_{|z|\leq\delta}\int_{\mR^d}\cR_j p^{\psi}_{s,t}(-x)\sD^{\phi(\cdot,z)}_1 u(s,x)\dif x\nu(\dif z).
\end{align*}
By the assumption \eqref{GH13}, one can choose $\delta$ small enough so that
$$
\|\nabla^\ell\phi(s,\cdot,z)\|_\infty\leq C|z|\leq\tfrac{1}{2},\ |z|<\delta,\ \ \ell=0,1.
$$ 
In particular, the assumption \eqref{SG0} is satisfied.
Let $\theta\in(\alpha-1,1)$. By \eqref{KF8} and H\"older's inequality, we have for any $j\in\mN$,
\begin{align*}
\int^t_0|\cR_jP^{\psi}_{s,t}\sA_\delta u(s,0)|\dif s
&\leq\int^t_0\int_{|z|\leq\delta}
\left|\int_{\mR^d}\cR_j p^{\psi}_{s,t}(-x)\sD^{\phi(\cdot,z)}_1 u(s,x)\dif x\right|\nu(\dif z)\dif s\\
&\lesssim\int^t_0\|u(s)\|_{\bC^\theta}
\Big[\mu_{1+\theta}\big(|\nabla^2\cR_j p^{\psi}_{s,t}|\big)^\theta\mu_{1+\theta}\big(|\nabla \cR_j p^{\psi}_{s,t}|\big)^{1-\theta}\\
&\qquad\qquad+\mu_0\big(|\cR_j p^{\psi}_{s,t}|\big)+\mu_1\big(|\nabla \cR_j p^{\psi}_{s,t}|\big)\Big]\dif s\\
&\lesssim\left(\int^t_0\mu_{1+\theta}\big(|\nabla^2\cR_j p^{\psi}_{s,t}|\big)\|u(s)\|_{\bC^\theta}\dif s\right)^\theta\\
&\quad\times\left(\int^t_0\mu_{1+\theta}\big(|\nabla\cR_j p^{\psi}_{s,t}|\big)\|u(s)\|_{\bC^\theta}\dif s\right)^{1-\theta}\\
&\quad+\int^t_0\Big[\mu_0\big(|\cR_j p^{\psi}_{s,t}|\big)+\mu_1\big(|\nabla\cR_j p^{\psi}_{s,t}|\big)\Big]\|u(s)\|_{\bC^\theta}\dif s\\
&=:I_1(t)+I_2(t).
\end{align*}
Let $\beta\in(0,\alpha)$. Since $\beta<\alpha<1+\theta$ and $\gamma+\eps<\alpha$, and recalling
$$
\mu_{1+\theta}(\dif x)=(|x|\wedge 1)^{1+\theta}\dif x\leq |x|^\beta\dif x= m_\beta(\dif x),
$$
for $j\in\mN$, by \eqref{KF1} and \eqref{LJ1}, we have
\begin{align*}
I_1(t)&\lesssim \left(\int^t_0 m_\beta\big(|\nabla^2\cR_j p^{\psi}_{s,t}|\big)\|u(s)\|_{\bC^\theta}\dif s\right)^\theta
\left(\int^t_0 m_\beta\big(|\nabla\cR_j p^{\psi}_{s,t}|\big)\|u(s)\|_{\bC^\theta}\dif s\right)^{1-\theta}\\
&\lesssim \left(2^{(2-\beta-\gamma-\eps) j}\int^t_0(t-s)^{-\frac{\gamma+\eps}{\alpha}}\|u(s)\|_{\bC^\theta}\dif s\right)^\theta\\
&\times\left(2^{(1-\beta-\gamma-\eps) j}\int^t_0(t-s)^{-\frac{\gamma+\eps}{\alpha}}\|u(s)\|_{\bC^\theta}\dif s\right)^{1-\theta}\\
&\lesssim2^{(1+\theta-\beta-\gamma-\eps) j}\int^t_0(t-s)^{-\frac{\gamma+\eps}{\alpha}}\|u(s)\|_{\bB^\theta_{\infty,\infty}}\dif s,
\end{align*}
and also,
$$
I_2(t)\lesssim 2^{-\gamma j}\int^t_0(t-s)^{-\frac{\gamma}{\alpha}}\|u(s)\|_{\bB^\theta_{\infty,\infty}}\dif s.
$$
For $j=0$, by \eqref{KF10}, we clearly have
$$
\int^t_0|\cR_0P^\psi_{s,t}\sA_\delta u(s,0)|\dif s\lesssim 
\int^t_0\|u(s)\|_{\bB^\theta_{\infty,\infty}}\dif s.
$$
Thus, we obtain that for any $j\in\mN_0$,
$$
\int^t_0|\cR_jP^\psi_{s,t}\sA_\delta u(s,0)|\dif s\lesssim 2^{(1+\theta-\beta-\eps-\gamma) j}
\int^t_0(t-s)^{-\frac{\gamma+\eps}{\alpha}}\|u(s)\|_{\bB^\theta_{\infty,\infty}}\dif s.
$$
In particular, if we choose $\theta$ close to $\alpha-1$ from above and $\beta$ close to $\alpha$ from below so that
$$
1+\theta-\beta\leq\eps, \ \ \theta-(\alpha-1)\leq \eps,
$$
then we get for any $j\in\mN_0$,
\begin{align}\label{DD0}
2^{\gamma j}\int^t_0|\cR_jP^\psi_{s,t}\sA_\delta u(s,0)|\dif s\lesssim \int^t_0(t-s)^{-\frac{\gamma+\eps}{\alpha}}\|u(s)\|_{\bB^{\alpha-1+\eps}_{\infty,\infty}}\dif s.
\end{align}
 Recall $m_\beta (\dif x)=|x|^\beta \dif x$.
For $\bar\sA_\delta u$, letting $\bar\phi(x,z):=\phi(x,z)-\phi(0,z)$, by Fubini's theorem and the integration by parts, we have 
\begin{align*}
&|\cR_jP^{\psi}_{s,t}\bar\sA_\delta u(0)|=\left|\int_{\mR^d} \cR_j p^{\psi}_{s,t}(x)\bar\sA_\delta u(x)\dif x\right|
\leq 2\|u\|_\infty\int_{|z|>\delta}\int_{\mR^d}|\cR_j p^{\psi}_{s,t}(x)|\dif x\nu(\dif z)\\
&\qquad+\left|\int_{|z|>\delta}\int_{\mR^d}\Big(\cR_j p^{\psi}_{s,t}(x) \div_x \bar\phi(x,z)+\bar\phi(x,z)\cdot\nabla\cR_j p^{\psi}_{s,t}(x)\Big)u(x)\dif x\nu(\dif z)\right|\\
&\quad\lesssim \|u\|_\infty m_0\Big(|\cR_j p^{\psi}_{s,t}|\Big)+\|u\|_\infty m_1\Big(|\nabla\cR_j p^{\psi}_{s,t}|\Big)\int_{|z|>\delta}|z|\nu(\dif z),
\end{align*}
where we have used that $|\div_x \bar\phi(x,z)|\lesssim|z|$ and $|\bar\phi(x,z)|\lesssim |x|\cdot |z|$.
Thus by \eqref{KF1}, we get for any $\gamma\in(0,\alpha)$,
$$
\int^t_0|\cR_jP^{\psi}_{s,t}\bar\sA_\delta u(s,0)|\dif s\leq C2^{-\gamma j}\int^t_0 (t-s)^{-\frac{\gamma}{\alpha}}\|u(s)\|_\infty\dif s,
$$
which together with \eqref{DD0} yields the desired estimate.
\end{proof}
In order to obtain the gradient estimate, we need the following lemma to improve the regularity to higher order.
\bl\label{Le45}
For any $T>0$, $\gamma\in[0,\alpha+\alpha\wedge 1)$ and $\theta\in(\alpha,2)$,
there is a constant $C=C(T,\Theta,\gamma,\theta)>0$ such that for all $j\in\mN_0$ and $t\in(0,T]$,
$$
2^{\gamma j}\int^t_0|\cR_jP^{\psi}_{s,t}\sA u(s,0)|\dif s\leq 
C\int^t_0 (t-s)^{-\frac{\gamma'}{\alpha}}\|u(s)\|_{\bB^\theta_{\infty,\infty}}\dif s,
$$
where $\gamma':=\frac{\gamma}{2}\1_{\gamma\leq 2}+(\gamma-1)\1_{\gamma>2}<\alpha$.
\el
\begin{proof}
Still we only consider $\alpha\in(1,2)$. Let $\theta\in(\alpha,2)$. Noting that
\begin{align*}
&|u(x+\phi(x,z))-u(x+\phi(0,z))-(\phi(x,z)-\phi(0,z))\cdot\nabla u(x)|\\
&\quad\leq |\phi(x,z)-\phi(0,z)|\int^1_0|\nabla u(x+s\phi(x,z)+(1-s)\phi(0,z))-\nabla u(x)|\dif s\\
&\quad\stackrel{\eqref{GH13}}{\lesssim}|x||z|\Big((|z|^{\theta-1}\|\nabla u\|_{\bC^{\theta-1}})\wedge\|\nabla u\|_\infty\Big)
\lesssim |x|(|z|^{\theta}\wedge |z|)\|u\|_{\bC^{\theta}},
\end{align*}
by \eqref{LJ1}, we have
$$
|\sA u(x)|\lesssim 
|x|\cdot\|u\|_{\bC^{\theta}}\int_{\mR^d}(|z|^\theta\wedge|z|)\nu(\dif z)\lesssim |x|\cdot\|u\|_{\bB^{\theta}_{\infty,\infty}}.
$$
Hence,
\begin{align*}
\int^t_0|\cR_jP^{\psi}_{s,t}\sA u(s,0)|\dif s
&\leq\int^t_0\!\!\int_{\mR^d} |\cR_j p^{\psi}_{s,t}(-x)\sA u(s,x)|\dif x\dif s\\
&\lesssim\int^t_0
m_1 (|\cR_j p^{\psi}_{s,t}|)\|u(s)\|_{\bB^{\theta}_{\infty,\infty}}\dif s,
\end{align*}
which implies by \eqref{KF1} that for any $\bar\gamma\in[0,\alpha)$,
$$
2^{(\bar\gamma+\bar\gamma\wedge 1) j}\int^t_0|\cR_jP^{\psi}_{s,t}\sA u(s,0)|\dif s\leq 
C\int^t_0 (t-s)^{-\frac{\bar\gamma}{\alpha}}\|u(s)\|_{\bB^\theta_{\infty,\infty}}\dif s.
$$
Here $m_1 (\dif x):= |x| \dif x$.
The proof is completed by the change of $\gamma=\bar\gamma+\bar\gamma\wedge 1$.
\end{proof}

Next comes to treat the drift term.
\bl\label{Le44}
Let $\alpha\in( {1}/{2}, 2)$, $\beta\in[0,\alpha\wedge 1)$ with $\alpha+\beta>1$. Under \eqref{BB},
for any $T>0$, $\gamma\in[0,\alpha)$ and $\eps\in(0,(\alpha-\gamma)\wedge\beta)$, 
there is a constant $C=C(T,\Theta,\gamma,\eps)>0$ such that for all $j\in\mN_0$ and $t\in(0,T]$,
\begin{align}\label{NH1}
2^{\gamma j}\int^t_0|\cR_jP^{\psi}_{s,t}(b\cdot\nabla u)(s,0)|\dif s\leq C\int^t_0 (t-s)^{-\frac{\gamma+\eps}{\alpha}}\|u(s)\|_{\bB^{1-\beta+\eps}_{\infty,\infty}}\dif s,
\end{align}
and for any $\gamma\in[0,\alpha+\beta)$,
\begin{align}\label{NH2}
2^{\gamma j}\int^t_0|\cR_jP^{\psi}_{s,t}(b\cdot\nabla u)(s,0)|\dif s\leq C\int^t_0 (t-s)^{-\frac{(\gamma-\beta)\vee 0}{\alpha}}\|\nabla u(s)\|_\infty\dif s.
\end{align}
\el

\begin{proof}
For $j\in\mN_0$, by definition and \eqref{TR1}, we have
\begin{align*}
&\cR_jP^{\psi}_{s,t}(b\cdot\nabla u)(0)=\!\!\int_{\mR^d}\!\!\cR_j p^{\psi}_{s,t}(-x)(b\cdot\nabla u)(x)\dif x=\!\!\int_{\mR^d}\!\!\bar\cR_j p^{\psi}_{s,t}(-x)\cR_j(b\cdot\nabla u)(x)\dif x\\
&\qquad=\int_{\mR^d}\bar\cR_j p^{\psi}_{s,t}(-x)(b\cdot\nabla\cR_j u)(x)\dif x+\int_{\mR^d}\bar\cR_j p^{\psi}_{s,t}(-x)[\cR_j,b\cdot\nabla]u(x)\dif x.
\end{align*}
By \eqref{KF9} and \eqref{GS1}, we have for any $\theta\in(1-\beta,1)$,
\begin{align*}
|\cR_jP^{\psi}_{s,t}(b\cdot\nabla u)(0)|&\lesssim\|\nabla\cR_j u\|_\infty m_\beta(|\bar\cR_j p^{\psi}_{s,t}|)
+\|[\cR_j,b\cdot\nabla]u\|_\infty\|\bar\cR_j p^{\psi}_{s,t}\|_1\\
&\lesssim 2^{(1-\theta)j}\|u\|_{\bC^{\theta}} m_\beta(|\bar\cR_j p^{\psi}_{s,t}|)
+2^{(1-\theta-\beta)j}\|u\|_{\bC^{\theta}}\|\bar\cR_j p^{\psi}_{s,t}\|_1,
\end{align*}
where $ m_\beta(\dif x)=|x|^\beta\dif x$.
Since $\gamma+\eps<\alpha$, by  \eqref{KF1} and \eqref{LJ1}, we obtain that for any $\theta\in(1-\beta,1)$,
$$
\int^t_0|\cR_jP^{\psi}_{s,t}(b\cdot\nabla u)(s,0)|\dif s\leq C2^{(1-\theta-\beta-\gamma-\eps)j}\int^t_0 (t-s)^{-\frac{\gamma+\eps}{\alpha}}\|u(s)\|_{\bB^{\theta}_{\infty,\infty}}\dif s,
$$
which implies \eqref{NH1} by choosing $\theta=1-\beta+\eps$.
Moreover,  we also have
\begin{align*}
&\int^t_0|\cR_jP^{\psi}_{s,t}(b\cdot\nabla u)(s,0)|\dif s\leq\int^t_0\!\!\int_{\mR^d}|\cR_jp^{\psi}_{s,t}(-x)(b\cdot\nabla u)(s,x)|\dif x\dif s\\
&\quad\stackrel{\eqref{KF9}}{\lesssim}\int^t_0 m_\beta(|\cR_j p^{\psi}_{s,t}|)\|\nabla u(s)\|_\infty\dif s
\stackrel{\eqref{KF1}}{\lesssim}2^{-(\bar\gamma+\beta) j}\int^t_0 (t-s)^{-\frac{\bar\gamma}{\alpha}}\|\nabla u(s)\|_\infty\dif s,
\end{align*}
where $\bar\gamma\in[0,\alpha)$.
Thus we obtain \eqref{NH2} by letting $\gamma=\bar\gamma+\beta$.
\end{proof}

\br\label{Br47}\rm 
The reason of requiring $\alpha\in({1}/{2},2)$ in Lemma \ref{Le44} is due to   
$\beta<\alpha\wedge1$ and $\alpha+\beta>1$.
Here we require $\beta<\alpha$
because of the moment estimate in \eqref{He1}.  
It is an open problem whether this restrict $\alpha>1/2$ can be dropped.
\er

Now we are in a position to give
\begin{proof}[Proof of Theorem \ref{Th31}] 
We  divide the proof into two steps.
\medskip\\
{\it Step 1}. In this step we prove \eqref{DJ1} for $\gamma\in[0,\alpha)$.
Let
$$
\delta:=(\alpha-1)\vee (1-\beta),\ \  \eta\in(\delta-\alpha,\gamma].
$$
By \eqref{KL1} and Lemmas \ref{Le42}, \ref{Le43} and \eqref{NH1}, for any $\eps\in(0,(\alpha-\gamma)\wedge\beta)$, we have
\begin{align}\label{NH3}
\|u(t)\|_{\bB^\gamma_{\infty,\infty}}
&=\sup_{j\geq 0}2^{\gamma j}\|\cR_j u(t)\|_\infty\stackrel{\eqref{AQ1}}{=}\sup_{j\geq 0}2^{\gamma j}\sup_{y}|\cR_j  u^y(t,0)|\no\\
&\lesssim  t^{\frac{\eta-\gamma}{\alpha}}\|\varphi\|_{\bB^\eta_{\infty,\infty}}
+\int^t_0 (t-s)^{-\frac{\gamma+\eps}{\alpha}}\|u(s)\|_{\bB^{\delta+\eps}_{\infty,\infty}}\dif s.
\end{align}
Since $\delta<\alpha$ and $\eta>\delta-\alpha$, one can choose $\eps\in(0,(\alpha-\gamma)\wedge\beta)$ small enough and $\eta'$ so that
$$
\delta+2\eps<\alpha,\ \ \eta'\in(\delta+\eps-\alpha,(\delta+\eps)\wedge\eta].
$$
Thus by \eqref{NH3} with $\gamma=\delta+\eps$ and $\eta=\eta'$, we have
$$
\|u(t)\|_{\bB^{\delta+\eps}_{\infty,\infty}}\lesssim  t^{\frac{\eta'-\delta-\eps}{\alpha}}\|\varphi\|_{\bB^{\eta'}_{\infty,\infty}}
+\int^t_0 (t-s)^{-\frac{\delta+2\eps}{\alpha}}\|u(s)\|_{\bB^{\delta+\eps}_{\infty,\infty}}\dif s,
$$ 
which implies by Gronwall's inequality (see Lemma \ref{Le26}) that  for all $t\in(0,T]$,
\begin{align}\label{KP2}
\|u(t)\|_{\bB^{\delta+\eps}_{\infty,\infty}}\lesssim t^{\frac{\eta'-\delta-\eps}{\alpha}}\|\varphi\|_{\bB^{\eta'}_{\infty,\infty}}
\lesssim t^{\frac{\eta'-\delta-\eps}{\alpha}}\|\varphi\|_{\bB^{\eta}_{\infty,\infty}}.
\end{align}
Now substituting \eqref{KP2} into \eqref{NH3}, we obtain that for all $t\in(0,T]$,
\begin{align}
&\|u(t)\|_{\bB^\gamma_{\infty,\infty}}\lesssim \|\varphi\|_{\bB^\eta_{\infty,\infty}}\left(t^{\frac{\eta-\gamma}{\alpha}}
+\int^t_0 (t-s)^{-\frac{\gamma+\eps}{\alpha}}s^{\frac{\eta'-\delta-\eps}{\alpha}}\dif s\right)\no\\
&\quad\lesssim  t^{\frac{\eta-\gamma}{\alpha}}\|\varphi\|_{\bB^\eta_{\infty,\infty}}
\Rightarrow
\|Q_{s,t}\varphi\|_{\bB^{\gamma}_{\infty,\infty}}\leq C (t-s)^{\frac{\eta-\gamma}{\alpha}}\|\varphi\|_{\bB^\eta_{\infty,\infty}}.\label{PW2}
\end{align}
Thus we obtain \eqref{DJ1} for any $\gamma\in[0,\alpha)$ since $\delta-\alpha=-((\alpha+\beta-1)\wedge 1)$.
\medskip
\\
{\it Step 2}. In this step we improve the spatial regularity of $Q_{s,t}\varphi$ to  order $\gamma\in[0,\alpha+\alpha\wedge\beta)$.
We consider the case of 
$\alpha>1/2$. The case of $\alpha\in(0,1/2]$ and $b\equiv 0$ is easier. 
Let 
$$
\gamma\in(\alpha\vee 1,\alpha+\alpha\wedge\beta),\  \  \eta\leq\gamma,\ \ \theta\in(\alpha\vee 1,2).
$$
 By Lemma 
\ref{Le45} and \eqref{NH2}, 
for  $\gamma'$ being as in Lemma \ref{Le45}, we have
\begin{align}\label{PW3}
\begin{split}
\|u(t)\|_{\bB^{\gamma}_{\infty,\infty}} \lesssim  & \, 
t^{\frac{\eta-\gamma}{\alpha}}\|\varphi\|_{\bB^\eta_{\infty,\infty}}
+\int^t_0 (t-s)^{-\frac{\gamma'}{\alpha}}\|u(s)\|_{\bB^\theta_{\infty,\infty}}\dif s\\
&+\int^t_0 (t-s)^{-\frac{\gamma-\beta}{\alpha}}\|\nabla u(s)\|_\infty\dif s.
\end{split}
\end{align}
In particular, for  $\gamma'':=\gamma'\vee(\gamma-\beta)<\alpha$, we have
$$
\|u(t)\|_{\bB^{\gamma}_{\infty,\infty}}\lesssim t^{\frac{\eta-\gamma}{\alpha}}\|\varphi\|_{\bB^\eta_{\infty,\infty}}
+\int^t_0 (t-s)^{-\frac{\gamma''}{\alpha}}\|u(s)\|_{\bB^\gamma_{\infty,\infty}}\dif s,
$$
which implies by Gronwall's inequality that for any $\eta\in(\gamma-\alpha,\gamma]$,
\begin{align}\label{PW1}
\|u(t)\|_{\bB^{\gamma}_{\infty,\infty}}\lesssim t^{\frac{\eta-\gamma}{\alpha}}\|\varphi\|_{\bB^\eta_{\infty,\infty}}\Rightarrow
\|Q_{s,t}\varphi\|_{\bB^{\gamma}_{\infty,\infty}}\lesssim (t-s)^{\frac{\eta-\gamma}{\alpha}}\|\varphi\|_{\bB^\eta_{\infty,\infty}}.
\end{align}
For $\eta\in(\delta-\alpha,\alpha)$, by \eqref{GQ1} and the estimate obtained in Step 1, we have
\begin{align*}
\|Q_{s,t}\varphi\|_{\bB^{\gamma}_{\infty,\infty}}=\|Q_{\frac{s+t}{2},t}Q_{s,\frac{s+t}{2}}\varphi\|_{\bB^{\gamma}_{\infty,\infty}}
&\stackrel{\eqref{PW1}}{\lesssim} (t-s)^{\frac{\eta'-\gamma}{\alpha}}\|Q_{s,\frac{s+t}{2}}\varphi\|_{\bB^{\eta'}_{\infty,\infty}}\\
&\stackrel{\eqref{PW2}}{\lesssim} (t-s)^{\frac{\eta-\gamma}{\alpha}}\|\varphi\|_{\bB^\eta_{\infty,\infty}},
\end{align*}
where $\eta'\in(\gamma-\alpha,\alpha)\subset(\gamma-\alpha,\gamma)$.
The proof is completed by interpolation.
\end{proof}

\section{Proof of Theorem \ref{Main}}

Now we give the proof of Theorem \ref{Main} under condition (ii). 
 Case (i) is easier.
 Thus we assume $\alpha\in( {1}/{2},2)$.
We divide the proof into three steps.
\medskip
\\
{\it {Step 1}}.  Fix $T>0$ and $\gamma\in(\alpha\vee 1,\alpha+\alpha\wedge 1)$. For any $\varphi\in \bC^{\gamma-}$, we first
show the existence of a classical solution $u^T\in C([0,T]; \bC^{\gamma-})$ for the following backward nonlocal-PDE 
\begin{align}\label{PDE11}
\p_su^T+\sL^{\sigma,b}_s u^T=0,\ u^T(T)=\varphi,
\end{align}
where $\sL^{\sigma,b}_s$ is defined by \eqref{SL}. For $n\in\mN$, define
$$
\varphi_n(x):=\varphi*\rho_n(x),\ b_n(t,x):=b(t,\cdot)*\rho_n(x),\ \ \sigma_n(t,x):=\sigma(t,\cdot)*\rho_n(x),
$$
 where $(\rho_n)_{n\in\mN}$ is a family of mollifiers in $\mR^d$.
 It is easy to see that 
$$
\varphi_n\in C^\infty_b(\mR^d),\ \ \sigma_n,\ b_n\in L^\infty(\mR_+; C^\infty_b(\mR^d)).
$$
It is well known that under these assumptions, for any $s\geq 0$ and $x\in\mR^d$, the following SDE admits a unique strong solution $X^n_{s,t}(x)$:
$$
\dif X^n_{s,t}= b_n(t,X^n_{s,t})\dif t+\sigma_n(t,X^n_{s,t-})\dif  Z_t,\ \ X^n_{s,s}=x.
$$
Moreover, $\{X^n_{s,t}(x),x\in\mR^d, t\geq s\geq 0\}$ forms a $C^\infty$-stochastic flow, and
$$
u^T_n(s,x):=\mE\varphi_n(X^n_{s,T}(x))\in C([0,T]; C^2_b(\mR^d))
$$
uniquely solves the following equation:
\begin{align}\label{HW1}
\p_su^T_n+\sL^{\sigma_n,b_n}_s u^T_n=0,\ u^T_n(T)=\varphi_n.
\end{align}
Below we let
$$
\phi_n(t,x,z):=\sigma_n(t,x) z,\ \ \nu(\dif z):=\sum_{i=1}^d |z_i|^{-1-\alpha}\delta_{0}(\dif z_1)\cdots\dif z_i\cdots\delta_{0}(\dif z_d).
$$
Under ({\bf H$^\sigma$}) and $\|\nabla\sigma\|_\infty<\infty$, it is easy to see that \eqref{GH13}-\eqref{BB} hold uniformly for the above $b_n,\phi_n$ and $\nu$.
Thus for any $\gamma'\in(\alpha\vee 1, \gamma)$, one can use \eqref{DJ1} to derive the following uniform estimate:
$$
\sup_n\|u^T_n\|_{C([0,T]; \bC^{\gamma'})}\leq C\|\varphi\|_{\bC^{\gamma'}}.
$$
By \eqref{HW1} and the above uniform estimate, one sees that for all $0\leq t_0<t_1\leq T$,
$$
\|u^T_n(t_1)-u^T_n(t_0)\|_\infty\leq \int^{t_1}_{t_0}\|\sL^{\sigma_n,b_n}_s u^T_n(s)\|_\infty\dif s\leq C|t_1-t_0|,
$$
where $C>0$ is independent of $n$.
Now by Ascolli-Arzela's lemma, there are function $u\in C([0,T]; \bC^{\gamma-})$ and subsequence still denoted by $n$ 
such that for any $T, R>0$,
$$
\lim_{n\to\infty}\|\nabla^j u^T_n-\nabla^ju^T\|_{C([0,T]\times B_R)}=0,\ \ j=0,1.
$$
Taking $n\to \infty$ in \eqref{HW1}, 
one finds that $u$ is a classical solution of nonlocal-PDE \eqref{PDE11}
in the sense of Definition \ref{Def41}.
\medskip\\
\textit {Step 2}.  
Let $u^T\in C([0,T]; \bC^{\gamma-})$ be the classical solution of nonlocal equation
\eqref{PDE11}. Let $X_{s,t}(x)$ be the unique solution of SDE \eqref{SDE0}.
By applying It\^o's formula to $(t,x)\mapsto u^T(t,x)$, we obtain
\begin{align*}
u^T(T,X_{s,T}(x))=u^T(s, x)+\int^T_s(\p_tu^T+\sL^{\sigma,b}_tu^T)(t,X_{s,t}(x))\dif t+\mbox{ a martingale}.
\end{align*}
Hence, by \eqref{PDE11},
$$
P_{s,T}\varphi(x)=\mE\varphi(X_{s,T}(x))=\mE u^T(T,X_{s,T}(x))=u^T(s,x).
$$
The desired estimate \eqref{DW1} now follows by Theorem \ref{Th31}.
\medskip\\
\textit {Step 3}.  
For  {\bf (A)}, let $\varphi\in \bB^{-\eta}_{\infty,\infty}$ for some $\eta<(\alpha+\beta-1)\wedge 1$ and let
$\varphi_\eps:=\varphi*\rho_\eps$ be the mollifying approximation. Clearly, by \eqref{PDE11} we have
$$
\p_s P_{s,t}\varphi_\eps(x)+\sL^{\sigma,b}_s P_{s,t}\varphi_{\eps}(x)=0.
$$
In particular, for any $0\leq t_0<t_1<t$ and $x\in\mR^d$,
$$
P_{t_0,t}\varphi_\eps(x)=P_{t_1,t}\varphi_\eps(x)+\int^{t_1}_{t_0}\sL^{\sigma,b}_s P_{s,t}\varphi_{\eps}(x)\dif s.
$$
By \eqref{DW1} and taking limits $\eps\to 0$, we obtain \eqref{DW01}.

For {\bf (B)}, since $\alpha>1/2$ and $\alpha+\beta>1$, 
one can choose $\gamma>1$ and $\eta=0$ in \eqref{DW1} so that
$$
\|P_{s,t}\varphi\|_{\bC^{\gamma}}\leq C(t-s)^{-\frac{\gamma}{\alpha}}\|\varphi\|_{\bB^0_{\infty,\infty}}\leq C(t-s)^{-\frac{\gamma}{\alpha}}\|\varphi\|_{\infty}.
$$
On the other hand, it is clear that
$$
\|P_{s,t}\varphi\|_{\infty}\leq \|\varphi\|_{\infty}.
$$
The desired gradient estimate now follows by interpolation theorem (see \cite[p35, Theorem 3.2.1]{Kr}). 

For {\bf (C)}, let  $\eta\in(\delta-\alpha,0)$. By \eqref{DW1}, we have
$$
\|P_{s,t}\varphi\|_\infty\leq C_{s,t}\|\varphi\|_{\bB^\eta_{\infty,\infty}}.
$$
From this, by Sobolev's embedding, one sees that 
$$
P_{s,t}\varphi(x)=\int_{\mR^d}\varphi(y)p_{s,t}(x,y)\dif y,\ \ p_{s,t}(x,\cdot)\in \bB^{-\eta}_{1,1}.
$$
Thus, we obtain the desired regularity.

\end{document}